\documentclass[12pt]{article}

\usepackage{amsmath}
\newcommand{\cpv}{\ensuremath{\!\!\!\!\!\!\!\!\! - \ \ \ }}
\newcommand{\cpvBox}{\ensuremath{\!\!\!\!\!\!\!\!\!\! - \ \ \ }}

\newcommand{\hypsng}{\ensuremath{\!\!\!\!\!\!\!\!\!\! = \ \ \ }}
\newcommand{\hpsngAbs}{\ensuremath{\!\!\!\!\!\!\!\!\!\!\! = \ \ \ }}
\newcommand{\hpIV}{\ensuremath{\!\!\!\!\!\!\!\!\!\! =  \ \ }}

\newcommand{\cpvtex}{\ensuremath{\!\!\!\!\!\! - \  }}
\newcommand{\hypsngtex}{\ensuremath{\!\!\!\!\!\! = \  }}
\newcommand{\hpTexI}{\ensuremath{\!\!\!\!\!\!\! = \  }}

\newcommand{\myLp}{\ensuremath{\ell^{\prime}}}
\newcommand{\myE}{\ensuremath{\textrm{\large{e}}}}

\title{\vspace*{-2.5cm}
       {\small \sf
         Submitted to ``Quarterly of Applied Mathematics''}
       \\
       \vspace*{1cm}{\bf \large{INTEGRAL EQUATIONS WITH \\
HYPERSINGULAR KERNELS -- THEORY AND APPLICATIONS TO FRACTURE MECHANICS}}}

\author{Youn-Sha Chan\footnote{Department of Mathematics,
University of California, Davis, CA 95616, U.S.A.}$^{\
,\ddag}$, Albert C. Fannjiang$^{\thefootnote ,\ddag}$, and
Glaucio H. Paulino\footnote{Department of Civil and Environmental
Engineering, University of Illinois, 2209 Newmark Laboratory, 205
North Mathews Avenue, Urbana, IL 61801, U.S.A.}$^{\
,}$\footnote{Graduate Group in Applied Mathematics (GGAM), University of
California, Davis, CA 95616, U.S.A.}}

\begin{document}

\input{epsf}
\input{psfig}

\maketitle

\begin{abstract}
Hypersingular integrals of the type
\[ I_{\alpha}(T_n,m,r) = 
\int_{-1}^{1}\hpsngAbs\frac{T_n(s)(1-s^2)^{m-\frac{1}{2}}}{(s-r)^\alpha}ds
\ , \ \ \ |r|<1  \]
and
\[ I_{\alpha}(U_n,m,r) = 
 \int_{-1}^{1}\hpsngAbs\frac{U_n(s)(1-s^2)^{m-\frac{1}{2}}}{(s-r)^\alpha}ds
\ , \ \ \ |r|<1  \]
are investigated for general integers $\alpha$ (positive)
and $m$ (non-negative),
where $T_n(s)$ and $U_n(s)$ are the Tchebyshev polynomials of the 1st
and 2nd kinds, respectively.  Exact formulas are derived for the cases
$\alpha = 1, 2, 3, 4$ and $m = 0, 1, 2, 3$; most of them
corresponding to new solutions derived in this paper.  Moreover, a
systematic approach
for evaluating these integrals when $\alpha > 4$ and $m>3$ is
provided.  The integrals are also evaluated as
$|r|>1$ in order to calculate stress intensity factors (SIFs).
Examples involving crack problems are given and discussed with
emphasis on the linkage between mathematics and mechanics of
fracture.  The examples include classical linear
elastic fracture mechanics (LEFM), functionally graded materials
(FGM), and gradient elasticity theory.  An appendix, with closed form
solutions for a broad class of integrals, supplements the paper.
\end{abstract}

\section{Introduction}
Finite and boundary element methods are two of the most
frequently used numerical approaches for solving crack problems in
fracture mechanics.  An alternative approach is the integral
equation method, which is more efficient ({\it{e.g.}} it reduces a
partial differential equation ({\sf PDE}) in two dimensions 
to an one-dimensional integral
equation) and, in general, is more accurate than the
aforementioned methods.
The accuracy of the integral equation method relies
on the analytical evaluation of singular kernels to cancel the
singularity (regularization).
In general, the cancellation of singularity is
not trivial, in particular, in the case of
hypersingular integrals.  This is the main concern of this paper.

Integral equations arising in static crack
problems in fracture mechanics are typically 
Fredholm integral equations of the form
\begin{equation}
\int_{c}^{d}\mbox{\sf{kernel}}(x,t)\ 
D(t)\ dt \ \ =\ \ p(x)\ , \ \quad c<x<d \ ,
\label{int_xt1}
\end{equation}
where {\sf kernel}$(x,t)$ is, in general, a singular function of
$(x,t)$; $D(t)$ is the unknown, called density function; $p(x)$ is
some known (input) function corresponding to the loading on the crack
faces; and the interval $(c, d)$ refers to the crack surfaces where
$2a = d - c$ denotes the crack length.
By the Fourier transform, we write
\begin{equation}
\mbox{\sf kernel}(x,t) \ \ =\ \
\int_{-\infty}^{\infty}K(\xi)\myE^{i(t-x)\xi}d\xi \ .
\end{equation}
The singular part of the kernel can be separated from the regular part,
by decomposing the Fourier transform as
\begin{equation}
K(\xi)=\underbrace{K_{\infty}(\xi)}_{\mbox{singular}} \ +\ 
\underbrace{[K(\xi) - K_{\infty}(\xi)]}_{\mbox{nonsingular}} \ ,
\label{decompose}
\end{equation}
which can be accomplished through asymptotic analysis (discussed later
in this paper).
Such an analysis is difficult for
complicated $K(\xi)$. This is another issue to be addressed  
in this paper.

Once the decomposition (\ref{decompose}) is
accomplished,
the integral equation (\ref{int_xt1}) can be rewritten as
\begin{equation}
\int_{c}^{d} \hpIV\frac{c_{\alpha}\ D(t)}{(t-x)^\alpha}\ dt \ +\ 
\int_{c}^{d}k(x,t)\ D(t)\ dt \ +\  f(x) \ = \ p(x)\ , \ \quad c<x<d \ ,
\label{int_xt2}
\end{equation}
where $\int\hpTexI$ denotes an improper integral; $c_{\alpha}$ is
a constant associated to the singular kernel
$1/(t-x)^\alpha$;
$k(x,t)$ is the nonsingular (regular) kernel; $f(x)$ is a function
standing for the free term; and $\alpha$ is a positive integer which
determines the degree of the singularity.
If $\alpha = 1$, the
integral equation (\ref{int_xt2}) is called a Cauchy singular
Fredholm integral equation, and the singular term is evaluated as
a Cauchy principal-value (CPV) integral.  If $\alpha \ge 2$,
it is called a hypersingular
Fredholm integral equation and
the singular term is evaluated as a Hadamard finite-part (HFP) integral
~\cite{KayaErdo87, Kutt75, Monegato94, Nedelec82}.
The notation $\int \cpvtex$ and $\int \hypsngtex$ refer to CPV and HFP, respectively.
Most of the works in the literature involve
either $\alpha = 1$~\cite{DelaErdo83,
Erdogan78, ErdoGupt72, ErdoGuptCook73, KondaErdo94, OztuErdo93}~ or
$\alpha = 2$~\cite{KayaDisrt84, KayaErdo87, SchuErdo98}.
Thus another focus of this paper is to deal with
hypersingular integral equations with $\alpha \geq 3$ which arise
naturally in gradient elasticity theories (see Example 3 in
\textsf{Section~\ref{exam}}).

The singular
and hypersingular integrals which involve Tchebyshev polynomials
($T_n$, first kind; $U_n$, second kind) and
weight function $(1-s^2)^{m-\frac{1}{2}}$, $m \geq 0$, with singularity
$\alpha \geq 1$ are defined by
\begin{equation}
I_{\alpha}(T_n,m,r) =
\int_{-1}^{1}\hypsng\frac{T_n(s)(1-s^2)^{m-\frac{1}{2}}}{(s-r)^\alpha}ds
\ , \ \ |r|<1 \ ,
\label{gnrlTn}
\end{equation}
and
\begin{equation}
I_{\alpha}(U_n,m,r) =
 \int_{-1}^{1}\hypsng\frac{U_n(s)(1-s^2)^{m-\frac{1}{2}}}{(s-r)^\alpha}ds
\ , \ \ |r|<1 \ .
\label{gnrlUn}
\end{equation}
The scope of this paper is as follows.
First, Cauchy singular integrals, {\it{i.e.}}
$\alpha = 1$, are evaluated and exact formulas are derived for general
$m$.  The new results here are closed form analytical solutions for
$I_1(T_n,m,r)$, $m \geq 1$ and $I_1(U_n,m,r)$, $m \geq 2$.
Once $I_1(T_n,m,r)$ and $I_1(U_n,m,r)$ are known, hypersingular
integrals $I_{\alpha}(T_n,m,r)$ and $I_{\alpha}(U_n,m,r)$,
$\alpha \geq 2$, can be found by successive differentiation (with
respect to $r$) in the sense of finite-part integrals; formulas for
$I_2(T_n,m,r)$, $m \geq 1$ and $I_2(U_n,m,r)$, $m \geq 2$ are
derived in this manner.  In the cases where $\alpha \geq 3$, evaluation
of hypersingular integrals becomes tedious and the formulas are
lengthy. Thus $I_{\alpha}(T_n,m,r)$ and $I_{\alpha}(U_n,m,r)$ are
provided only for $\alpha = 3, 4$ and general $m$.

\section{Related Work}
Singular integral equations have played an active role in the field of
solid mechanics, particularly in the solution of fracture mechanics
problems.  For instance, the Journal {\it{Integral Equations and
Operator Theory}} has dedicated a special issue to
``Integral Equations Methods in Engineering and Physics'' (Volume 5,
No. 4, 1982).  Also, in June 1984, IMACS (International Association
for Mathematics and Computers in Simulation) has held a
Symposium~\cite{IMACS84} devoted to
``Numerical Solution of Singular Integral Equations''.

\ 

According to the notation introduced in equation (\ref{int_xt2}),
singular integral equations can be classified by the order of
singularity $\alpha$.  The case $\alpha = 1$ has been widely used and
well developed~\cite{ErdoGupt72, ErdoGuptCook73}.  A rich field of
application of singular integral equations is fracture mechanics of
bimaterial and nonhomogeneous materials.  For instance,
the investigation of crack
behavior in nonhomogeneous materials has found many applications to
functionally graded materials (FGMs)~\cite{DelaErdo83, Erdogan85,
JinNoda94}.  Another use of singular integral equations involves FGMs
for high temperature applications, so that thermal stress intensity
factors can be numerically calculated~\cite{JinBatra96c, NodaJin93}.
Application of hypersingular Fredholm integral equations  for $\alpha
\geq 2$ can be found in references~\cite{KayaDisrt84, KayaErdo87, SchuErdo98}.

\ 

Quadrature formulas which
involve hypersingular integrals have been drawing a considerable amount of
concentration~\cite{HuiShia99, KayaDisrt84, KayaErdo87, Korsunsky98, Krenk75b,
Martin96, Monegato94, TheoIoak77} after Kutt first introduced the Hadamard
finite-part (HFP) idea in
his work~\cite{Kutt75}.  Based on some previous work,
Kaya~\cite{KayaDisrt84} has presented a very nice
interpretation about HFP integrals.  One of the key steps in the
derivations involves the fact that higher order
singular integrals can be obtained from lower order ones
by exchangeability of differentiation and
integration~\cite{KayaErdo87, Monegato94}, {\it e.g.}
\begin{equation}
\int_{-1}^{1}\hypsng\frac{D(s)}{(s-r)^{\alpha+1}} ds =
\frac{1}{\alpha}\int_{-1}^{1}\hypsng\frac{\partial}{\partial r}\left[\frac{
D(s)}{(s-r)^{\alpha}} \right] ds=
\frac{1}{\alpha}\frac{d}{dr}\int_{-1}^{1}\hypsng\frac{D
(s)}{(s-r)^{\alpha}} ds
 \ , \quad |r|<1,
\label{diffKey}
\end{equation}
where $D(s)$ is the normalized density function.
For instance, in order to find
\[
\int_{-1}^{1}\hypsng\frac{D(s)}{(s-r)^2} ds \ , \quad |r|<1 \ ,
\]
it suffices to know how to evaluate
\[
\int_{-1}^{1}\cpv\frac{D(s)}{s-r} ds \ , \quad |r|<1 \ .
\]
This concept is applied later in this paper.

\

Another main motivation for numerical
evaluations of hypersingular integrals is
due to the boundary element method, and the reader is directed to the
review paper by Tanaka {\it et al}~\cite{TanakSlad94}.  Most
recent work has been focused
on singularity with $\alpha =2$ in two-dimensional problems.






\section{Theoretical Aspects}
First, relevant concepts involving integration and approximation are
given.  These concepts position the contribution of the work with
respect to the available literature.  Next, a discussion on the
influence of the density function on the corresponding singular
integral equation formulation is presented.  Afterwards,
basic properties of the
Tchebyshev polynomials are provided.  These properties are heavily
used in the analytical derivations that follow.


\subsection{Integration and Approximation}
As far as the integration and numerical procedures are concerned,
the integral equation (\ref{int_xt2}) may be normalized through
the following change of variables
\begin{equation}
s=\frac{2}{d-c}\left(t-\frac{c+d}{2}\right)\ \ \  \mbox{and}\ \ \
r=\frac{2}{d-c}\left(x-\frac{c+d}{2}\right)\ \ ,
\label{chngVar}
\end{equation}
which leads to the normalized version of the integral equation
(\ref{int_xt2})
written as\footnote{The notations in this paper have been chosen as
following: $x$ and $t$ refer to the physical quantities and have
dimension of ``Length''; $r$ and $s$ are normalized (dimensionless)
variables, corresponding to $x$ and $t$, respectively.}
\begin{equation}
\int_{-1}^{1}\hypsng\frac{D(s)}{(s-r)^\alpha}\ ds +
\int_{-1}^{1}\mathcal{K}(r,s)D(s)\ ds  + F(r) \ = \ P(r)\ , \ \
-1<r<1 \ .
\label{intNorm}
\end{equation}
The density function $D(s)$ is further assumed to have the representation
\begin{equation}
D(s) \ = \ R(s)W(s) \ .
\label{repDnsty}
\end{equation}
The weight function $W(s)$ determines the singular behavior of the
solution $D(s)$ and has the form
\begin{equation}
W(s)\ = \ (1-s)^{m_1}(1+s)^{m_2} \ \ .
\label{rep12Wt}
\end{equation}
In general, $m_1 \neq m_2$, and the corresponding integrals, which involve Jacobi
polynomials $P_n^{(m_1,m_2)}(s)$, are of the type
\begin{equation}
\int_{-1}^1\cpv\frac{(1-s)^{m_1}(1+s)^{m_2}P_n^{(m_1,m_2)}(s)}{s-r}ds
\ ,
\end{equation}
and can be expressed in terms of gamma and hypergeometric
functions~\cite{AbraSteg72, KayaDisrt84, KayaErdo87, Krenk75a}.  In this
paper, only the case $m_1 = m_2$ is considered and $m_1$, $m_2$ are
set to be
\begin{equation}
m_1 = m_2 = m-\frac{1}{2} \ .
\end{equation}
Thus $W(s)$ can be expressed as
\begin{equation}
W(s)\ = \ \left(1-s^2 \right)^{m-\frac{1}{2}} \ \
m = 0,1,2, \cdots \ \ .
\label{repWt}
\end{equation}

\ 

According to function-theoretic method~\cite{Erdogan78,
ErdoGuptCook73, KayaDisrt84, KayaErdo87, Muskhelishvili53}~, 
the value of $m$ is determined by the order of singularity $\alpha$.
For example, if $\alpha =1$, then $m =0$ and the fundamental solution
$D(s)$ to the Cauchy singular integral equation
(\ref{intNorm}) takes the form
\begin{equation}
D(s) \ = \ R(s)(1-s^2)^{-\frac{1}{2}}\ = \ \frac{R(s)}{\sqrt{1-s^2}}  \ .
\label{rep1Dnsty}
\end{equation}
In this case, which consists of the majority of the work involving
applications of 
integral equations to fracture mechanics~\cite{DelaErdo83,
Erdogan78, ErdoGupt72, ErdoGuptCook73, KondaErdo94, OztuErdo93}~,
$R(s)$ is chosen to be
\begin{equation}
R(s)\ = \ \sum_n^{\infty}a_nT_n(s)  \ \ ;
\label{expnd1Cheby}
\end{equation}
and because of that, the CPV integral $I_1(T_n,0,r)$ can be evaluated
exactly~\cite{AbraSteg72, KayaErdo87}~:
\begin{equation}
I_1(T_n,0,r) \ =\ 
\int_{-1}^{1}\cpv\frac{T_n(s)}{(s-r)\sqrt{1-s^2}}\ ds \
=\ \pi U_{n-1}(r) \ , \ \ \ \ n\geq 1\ .
\label{int1Cheby}
\end{equation}
Another reason for choosing the approximation (\ref{rep1Dnsty})
is that with respect to the weight function $W(s) = 1/ \sqrt{1-s^2}$, the
class of the Tchebyshev polynomials of first kind $T_n(s)$ is an orthogonal
family~\cite{AbraSteg72, KayaErdo87}:
\begin{equation}
\int_{-1}^1\frac{T_m(s)T_n(s)}{\sqrt{1-s^2}} ds = \left\{
   \begin{array}{lll}
\pi            &  m = n = 0   \\
\pi/2  &  m = n ; \ \ \ m,n = 1,2,3,  \cdots \\
0              &  m\neq n ; \ \ \ m,n = 0,1,2, \cdots
   \end{array}
\label{orthoCheby1}
                                         \right.
\end{equation}
With this orthogonal property a Galerkin-type method~\cite{Krenk75a}
may be applied to find the coefficients $a_n$ in equation (\ref{expnd1Cheby}).

\

If $\alpha =2$, then $m = 1$, and the solution $D(s)$ to the
hypersingular integral equation (\ref{intNorm}) is characterized by
\begin{equation}
D(s) \ = \ R(s)(1-s^2)^{\frac{1}{2}}\ = \ {R(s)}\sqrt{1-s^2}  \ .
\label{rep2Dnsty}
\end{equation}
Correspondingly, $R(s)$ is chosen to be
\begin{equation}
R(s)\ = \ \sum_n^{\infty}b_nU_n(s)  \ \ ,
\label{expnd2Cheby}
\end{equation}
because of the same reasons for the case $\alpha =1$, namely,
analytical evaluation and orthogonal property. With respect to the
first reason, the HFP
integral $I_2(U_n,1,r)$ can be evaluated analytically~\cite{KayaErdo87}~:
\begin{equation}
I_2(U_n,1,r) \ =\ 
\int_{-1}^{1}\hypsng\frac{U_n(s)\sqrt{1-s^2}}{(s-r)^2}ds =
-(n+1)\pi U_{n}(r)\ , \ \ \ \ n\geq 0 \ .
\label{int2UCheby}
\end{equation}
According to the second reason, by orthogonality,
\begin{equation}
\int_{-1}^1U_m(s)U_n(s)\sqrt{1-s^2} ds = \left\{
   \begin{array}{ll}
\pi/2  &  m = n ; \ \ \ m,n = 0,1,2, \cdots \\
0              &  m\neq n ; \ \ \ m,n = 0,1,2, \cdots\ \ ,
   \end{array}
\label{orthoCheby2}
                                         \right.
\end{equation}
and one may apply Galerkin-type methods~\cite{Krenk75a} to
find the coefficients $b_n$ in equation (\ref{expnd2Cheby}).

\ 

When $m = 3$, then $W(s) = (1-s^2)^{5/2}$, and neither $T_n(s)$ nor
$U_n(s)$ is an orthogonal family.  {\it However, if collocation method
is applied, one does not need the orthogonal property, as long as the
expansion function $R(s)$ is chosen such that}
\[
\int_{-1}^{1}\hypsng\frac{R(s)(1-s^2)^{\frac{5}{2}}}{(s-r)^\alpha}ds
\]
{\it can be evaluated analytically.}
For example, if $R(s)$ is expandeded as a Tchebyshev polynomial of the 1st
kind $T_n(s)$ or the 2nd kind $U_n(s)$, {\it{i.e.}}
\begin{equation}
R(s)\ = \ \sum_n^{\infty}a_nT_n(s) \ \ \ \mbox{or} \ \ \
R(s)\ = \ \sum_n^{\infty}b_nU_n(s) \ \ ,
\label{expndCheby}
\end{equation}
then the evaluation of
\[
I_{\alpha}(T_n,m,r) = \int_{-1}^{1}\hypsng\frac{T_n(s)\left(1-s^2
\right)^{m-\frac{1}{2}}}
{(s-r)^\alpha}\ ds \ \  \mbox{or} \ \ 
I_{\alpha}(U_n,m,r) = \int_{-1}^{1}\hypsng\frac{U_n(s)\left(1-s^2
\right)^{m-\frac{1}{2}}}
{(s-r)^\alpha}\ ds
\]
for general $m = 0,1,2,\cdots$ and $\alpha = 1,2,3,\cdots$ is a
necessary step for the numerical approach to the integral equation
(\ref{intNorm}).  This is the one of main tasks in this paper and is addressed
in {\textsf{Sections 4}} and {\textsf{5}} .



\subsection{Selection of the Density Function}
Usually the unknown function $D(t)$ in equation (\ref{int_xt1}) can be chosen as
the displacement profile ({\it{e.g.}} $u(t)$ -- a displacement function), the (first)
derivative of the displacement function ($du(t)/dt$, denoted by $\phi
(t)$ -- the slope function), or a higher derivative of $u(t)$.  The
choice of the unknown function $D(t)$ will affect the degree of
singularity in the formulation.  For example, consider the standard
mode III crack problem in a free space~\cite{Gdou90} and a 
linear elastic fracture mechanics (LEFM) setting.  If $D(t)$ is
chosen to be the slope function $\phi(t)$, then the governing integral
equation is the Cauchy
singular integral equation
\begin{equation}
D(t) \equiv \phi(t) \ \ , \ \ \ \int_{c}^{d}\cpv\frac{\phi(t)}{t-x}dt =
p(x)\ , \ \ c<x<d \ .  \label{CauchyForm}
\end{equation}
However, if $D(t)$ is chosen to be the displacement function $w(t)$, then the
hypersingular integral equation with $\alpha = 2$ is obtained,
\begin{equation}
D(t) \equiv w(t) \ \ , \ \ \ \int_{c}^{d}\hypsng\frac{w(t)}{(t-x)^2}dt =
p(x)\ , \ \  c<x<d \ . \label{QuadraForm}
\end{equation}
The differences between the above two formulations are discussed next.

\

In general, the (numerical) solution of a Cauchy singular integral
equation, {\it e.g.} equation (\ref{CauchyForm}), is easier than the
corresponding hypersingular integral equation, {\it e.g.} equation
(\ref{QuadraForm}).  A quick observation is that if equation
(\ref{CauchyForm}) is used, then the actual crack surface displacements are
obtained directly; while if equation (\ref{QuadraForm}) is used, an
extra step of integration is needed to recover $w(t)$.  However,
integration is not an unpleasant thing to do.

\

So, what can be gained from a more singular equation such as
(\ref{QuadraForm})? A formulation with more singular integral
may lead to a simpler kernel function and, thus,
simplify the
kernel evaluation and decomposition described in equation
(\ref{decompose}). Issues regarding the differences between
formulations of the type given in equations (\ref{CauchyForm}) and
(\ref{QuadraForm}) will be discussed in detail in the example section
of this paper.


\subsection{Properties of Tchebyshev Polynomials}
The evaluation of Cauchy singular and hypersingular integrals which
involve the Tchebyshev polynomials $T_n(s)$ and $U_n(s)$ highly depends on
the special properties of these polynomials.
They are listed here for the sake of completeness and because they
will be of much use later in the development of this work.
Most of them (but not all) can be
found in Hochstrasser~\cite{AbraSteg72} and Kaya and
Erdogan~\cite{KayaErdo87}.

\begin{itemize}
\item{Definition of Tchebyshev polynomials of the first
kind}:
\vspace{-.03in}
\begin{equation}
T_n(s) = \cos [n \cos^{-1}(s)] \ , \ \ \ \ n = 0, 1, 2, \cdots
\label{defCheby1}
\end{equation}

\item{Definition of Tchebyshev polynomials of the second
kind}:
\vspace{-.03in}
\begin{equation}
U_n(s) = \frac{\sin [(n+1)\cos^{-1}(s)]}{\sin [ \cos^{-1}(s)]}  \ , \
\ \ \ n = 0, 1, 2, \ldots
\label{defCheby2}
\end{equation}

\item{Iterative (recursive) properties}:
\vspace{-.06in}
\begin{equation}
sT_n(s) = \frac{1}{2}[T_{n+1}(s) + T_{n-1}(s)] \ , \ \ \ \ n\geq 1
\label{recur1Cheby}
\end{equation}

\begin{equation}
sU_n(s) = \frac{1}{2}[U_{n+1}(s) + U_{n-1}(s)] \ , \ \ \ \ n\geq 1
\label{recur2Cheby}
\end{equation}

\begin{equation}
T_n(s) = \frac{1}{2}[U_{n}(s) - U_{n-2}(s)] \ , \ \ \ \ n\geq 2
\label{recur3Cheby}
\end{equation}

\begin{equation}
U_{n}(s)(1-s^2) = sT_{n+1}(s) - T_{n+2}(s) \ , \ \ \ \ n\geq 0
\label{recur4Cheby}
\end{equation}

By means of equation (\ref{recur1Cheby}), one may rewrite equation
(\theequation) above as
\begin{equation}
U_{n}(s) = \frac{1}{2(1-s^2)}[T_{n}(s) - T_{n+2}(s)] \ , \ \ \ \ n\geq 0
\label{recur5Cheby}
\end{equation}

\ 

Thus an additional equality, which is useful in handling cubic
hypersingular integrals can be derived\footnote{The
equation number is stacked above the equal sign to show how the
equations are being derived and connected.}:
\begin{eqnarray}
\lefteqn{U_n(s)(1-s^2)^{\frac{3}{2}} 
 \stackrel{(\ref{recur5Cheby})}{=}
 \frac{1}{2}\left[T_{n}(s)-T_{n+2}(s)\right]\sqrt{1-s^2}} \nonumber \\
 & \stackrel{(\ref{recur3Cheby})}{=} & -\frac{1}{2}
 \left[\frac{1}{2}U_{n+2}(s)-U_n(s)+\frac{1}{2}U_{n-2}(s)\right]\sqrt{1-s^2}\ ,\
 \ n\geq 2  \nonumber \\
 & = & -\frac{1}{4}
 \left[U_{n+2}(s)-2U_n(s)+U_{n-2}(s)\right]\sqrt{1-s^2}\ ,\ \ n\geq 2 
\label{recur6Cheby}
\end{eqnarray}

\vspace{.5in}

\item{Derivatives}: \vspace{-.0in}
\begin{equation}
\frac{d T_n(s)}{ds} = nU_{n-1}(s) \ , \ \ \ \ n\geq 1
\label{derv1Cheby}
\end{equation}
\begin{equation}
\frac{d U_n(s)}{ds} = \frac{1}{1-s^2}\left[\frac{n+2}{2}U_{n-1}(s)
- \frac{n}{2}U_{n+1}(s)\right] \ , \ \ \ \ n\geq 1
\label{derv2Cheby}
\end{equation}

\end{itemize}


\section{Cauchy Singular Integral Formulas ($\alpha =1$)}
This section mainly evaluates $I_1(T_n,m,r)$ and $I_1(U_n,m,r)$, which
are defined in equations (\ref{gnrlTn}) and (\ref{gnrlUn}).  The
new result here is that the singular integral formulas are found for
general $m$.  In order to obtain this new result,
two well known Cauchy singular integral formulas are
introduced~\cite{AbraSteg72, KayaErdo87}: one is already stated in equation (\ref{int1Cheby}),
and the other one is
\begin{equation}
I_1(U_n,1,r)=\int_{-1}^{1}\cpv\frac{U_n(s)\sqrt{1-s^2}}{s-r}ds =
-\pi T_{n+1}(r)\ , \ \ \ \ n\geq 0 \ ,
\label{int2Cheby}
\end{equation}
which can be obtained as follows
\begin{eqnarray*}
\lefteqn{I_1(U_n,1,r)=\int_{-1}^{1}\cpv\frac{U_n(s)\sqrt{1-s^2}}{s-r}ds}\\
 & \stackrel{(\ref{recur5Cheby})}{=} &
 \frac{1}{2}\int_{-1}^{1}\cpv\frac{T_n(s)-T_{n+2}(s)}{\sqrt{1-s^2}(s-r)}ds \\
 & \stackrel{(\ref{int1Cheby})}{=} & \frac{\pi}{2}[U_{n-1}(r) -
 U_{n+1}(r)] \\
 & \stackrel{(\ref{recur3Cheby})}{=} & -\pi T_{n+1}(r)\ .
\end{eqnarray*}
The integral formulas for $m = 0,1,2,3, \cdots$ are derived below.
The general formulas have the restriction of minimum $n$.  The lower
$n$ terms can not be derived by general formulas, and are given in \textsf{Appendix A}.

\subsection{$I_1(T_n,m,r), \ m=0,1,2,3$}

\begin{itemize}
\item{$I_1(T_n,0,r)$}:\\
This is equation (\ref{int1Cheby}).
\item{$I_1(T_n,1,r)$}:
\begin{eqnarray}
\int_{-1}^{1}\cpv\frac{T_n(s)\sqrt{1-s^2}}{s-r}ds
 & \stackrel{(\ref{recur3Cheby})}{=} &
\frac{1}{2}\int_{-1}^{1}\cpv\frac{[U_n(s)-U_{n-2}(s)]\sqrt{1-s^2}}{s-r}ds
 \nonumber \\ 
 & \stackrel{(\ref{int2Cheby})}{=} & \frac{\pi}{2}[T_{n-1}(r)-T_{n+1}(r)]\ , \ \ \ \ n\geq 2 \ .
\label{intT11}
\end{eqnarray}

\item{$I_1(T_n,2,r)$}:
\begin{eqnarray}
\lefteqn{\int_{-1}^{1}\cpv\frac{T_n(s)(1-s^2)^{\frac{3}{2}}}{s-r}ds
 \stackrel{(\ref{recur3Cheby})}{=}
 \int_{-1}^{1}\cpv\frac{\frac{1}{2}[U_n(s)-U_{n-2}(s)](1-s^2)^{\frac{3}{2}}}{s-r}ds }
 \nonumber \\
 & \stackrel{(\ref{intU12})}{=} & \frac{\pi}{8}\left\{[T_{n-1}(r) - 2T_{n+1}(r) + T_{n+3}(r)] -
 [T_{n-3}(r) - 2T_{n-1}(r) + T_{n+1}(r)]\right\} \nonumber \\
 & = & -\frac{\pi}{8}[T_{n-3}(r) - 3T_{n-1}(r) + 3T_{n+1}(r) -
 T_{n+3}(r)]\ , \ \ \ \ n\geq 4
\label{intT12}
\end{eqnarray}

\item{$I_1(T_n,3,r)$}:
\begin{eqnarray}
\lefteqn{\int_{-1}^{1}\cpv\frac{T_n(s)(1-s^2)^{\frac{5}{2}}}{s-r}ds
 \stackrel{(\ref{recur3Cheby})}{=}
 \int_{-1}^{1}\cpv\frac{\frac{1}{2}[U_n(s)-U_{n-2}(s)](1-s^2)^{\frac{5}{2}}}{s-r}ds }
 \nonumber \\
 & \stackrel{(\ref{intU13})}{=} &
 \frac{\pi}{32}\left\{[T_{n-5}(r)-4T_{n-3}(r)+6T_{n-1}(r)-4T_{n+1}(r)+T_{n+3}(r)]\right.
 - \nonumber \\
 & \  & \qquad
 \left.[T_{n-3}(r)-4T_{n-1}(r)+6T_{n+1}(r)-4T_{n+3}(r)+T_{n+5}(r)]\right\}\ , \ \ n\geq 6
 \nonumber \\
 & = & \frac{\pi}{32}[T_{n-5}(r)-5T_{n-3}(r)+10T_{n-1}(r) \nonumber \\
 & \ & \ \hspace{.3in} -10T_{n+1}(r)+5T_{n+3}(r)-T_{n+5}(r)]\ , \ \ \ n\geq 6
\label{intT13}
\end{eqnarray}

\end{itemize}

\subsection{$I_1(U_n,m,r), \ m=1,2,3$}
\begin{itemize}
\item{$I_1(U_n,1,r)$}:\\
This is equation (\ref{int2Cheby}).

\item{$I_1(U_n,2,r)$}:
\begin{eqnarray}
\lefteqn{\int_{-1}^{1}\cpv\frac{U_n(s)(1-s^2)^{\frac{3}{2}}}{s-r}ds
  \stackrel{(\ref{recur5Cheby})}{=} 
\frac{1}{2}\int_{-1}^{1}\cpv\frac{[T_n(s)-T_{n+2}(s)]\sqrt{1-s^2}}{s-r}ds}
 \nonumber \\ 
 & \stackrel{(\ref{intT11})}{=} & \frac{\pi}{4}\left\{[T_{n-1}(r)-T_{n+1}(r)] -
 [T_{n+1}(r)-T_{n+3}(r)]\right\}\ , \quad n\geq 2 \nonumber \\
 & = & \frac{\pi}{4}[T_{n-1}(r) - 2T_{n+1}(r) + T_{n+3}(r)]\ , \quad n\geq 2 
\label{intU12}
\end{eqnarray}

\item{$I_1(U_n,3,r)$}:
\begin{eqnarray}
\lefteqn{\int_{-1}^{1}\cpv\frac{U_n(s)(1-s^2)^{\frac{5}{2}}}{s-r}ds
  \stackrel{(\ref{recur5Cheby})}{=} 
\frac{1}{2}\int_{-1}^{1}\cpv\frac{[T_n(s)-T_{n+2}(s)](1-s^2)^{\frac{3}{2}}}{s-r}ds}
 \nonumber \\ 
 & \stackrel{(\ref{intT12})}{=} &
  \frac{\pi}{16}\left\{[T_{n-1}(r)-3T_{n+1}(r)+3T_{n+3}(r)-T_{n+5}(r)]\right. - \nonumber \\
 & \ & \left.\qquad [T_{n-3}(r)-3T_{n-1}(r)+3T_{n+1}(r)-3T_{n+3}(r)]\right\}\
  ,\quad n\geq 4 \nonumber \\
 & = & -\frac{\pi}{16}[T_{n-3}(r)-4T_{n-1}(r)+6T_{n+1}(r)-4T_{n+3}(r)+T_{n+5}(r)]\ , \ \ n\geq 4
\label{intU13}
\end{eqnarray}

\end{itemize}

\subsection{$I_1(T_n,m,r)$ and $I_1(U_n,m,r)$}
At this point one may easily see the procedural steps above,
which take advantage of recursive properties
(\ref{recur3Cheby}) and (\ref{recur5Cheby}) between the Tchebyshev
polynomials $T_n(s)$ and $U_n(s)$. 
For instance, evaluation of
$I_1(T_n,4,r)=\int_{-1}^{1}\cpv T_n(s)(1-s^2)^{7/2}/(s-r)ds$
can be reduced to evaluation of
$I_1(U_n,4,r)=\int_{-1}^{1}\cpv U_n(s)(1-s^2)^{7/2}/(s-r)ds$,
which, in turn, can be reduced to evaluation of
$I_1(T_n,3,r)=\int_{-1}^{1}\cpv T_n(s)(1-s^2)^{5/2}/(s-r)ds$.
After a suitable number of steps, this reduction leads to
either (\ref{int1Cheby}) or (\ref{int2Cheby}).
This procedure is summarized in 
{\sf Figure~\ref{fig0}}.

Thus, by induction, one obtains the following formulas for $m\geq 1$.
\begin{figure}[ht]
{\Huge
\[
\boxed{
\begin{array}{c}
{\int_{-1}^{1}\cpvBox\frac{T_n(s)(1-s^2)^{m+\frac{1}{2}}}{s-r}ds}\\
\Downarrow \mbox{\normalsize{(\ref{recur3Cheby})}}\\
{\int_{-1}^{1}\cpvBox\frac{U_n(s)(1-s^2)^{m+\frac{1}{2}}}{s-r}ds}\\
\Downarrow \mbox{\normalsize{(\ref{recur5Cheby})}}\\
{\int_{-1}^{1}\cpvBox\frac{T_n(s)(1-s^2)^{m-\frac{1}{2}}}{s-r}ds}\\
\Downarrow \mbox{\normalsize{(\ref{recur3Cheby})}}\\
{\int_{-1}^{1}\cpvBox\frac{U_n(s)(1-s^2)^{m-\frac{1}{2}}}{s-r}ds}\\
\Downarrow \mbox{\normalsize{(\ref{recur5Cheby})}}\\
\hspace{-.25in} \mathbf{\vdots}\\
\mbox{\Large{Equation  (\ref{int1Cheby}) or (\ref{int2Cheby})}}
\end{array}
}
\]
}
\caption{Evaluation of $I_1(T_n,m,r)$ and $I_1(U_n,m,r)$ for
general $m$. The procedure reduces integrals with higher $m$
to those with lower $m$.}
 \label{fig0}
\end{figure}


\begin{itemize}
\item{$I_1(T_n,m,r)$, where $m\ge 1$, and $n\geq 2m$}
\begin{equation}
\boxed{
\int_{-1}^{1}\cpv\frac{T_n(s)(1-s^2)^{m-\frac{1}{2}}}{s-r}ds =\pi
(-1)^{m+1}\left(\frac{1}{2}\right)^{2m-1}\sum_{j=0}^{2m-1}(-1)^{j}\left(
\hspace{-.05in}\begin{array}{c} ^{2m-1}
\\ _{j} \end{array} \hspace{-.05in}\right)T_{n+1-2m+2j}(r)\ .}
\label{boxI}
\end{equation}

\item{$I_1(U_n,m,r)$, where $m\ge 2$, and $n\geq 2m-2$}
\begin{equation}
\boxed{
\int_{-1}^{1}\cpv\frac{U_n(s)(1-s^2)^{m-\frac{1}{2}}}{s-r}ds = \pi
(-1)^{m}\left(\frac{1}{2}\right)^{2m-2}\sum_{j=0}^{2m-2}(-1)^{j}\left(
\hspace{-.05in}\begin{array}{c} ^{2m-2}
\\ _{j} \end{array} \hspace{-.05in}\right)T_{n+3-2m+2j}(r)\ .}
\label{boxII}
\end{equation}
\end{itemize}
The usual notation
\[
\left(
\hspace{-.05in}\begin{array}{c} ^{m}
\\ _{j} \end{array} \hspace{-.05in}\right) =
\frac{(m)!}{j!\ (m-j)!}
\]
denotes the binomial coefficients.

\section{Hypersingular Integral Formulas ($\alpha \geq 2$)}
Once a Cauchy singular integral formula has been
reached, all other hypersingular integral formulas may be obtained
successively by
taking differentiation with respect to $r$, and making use of the
finite-part integral formula (\ref{diffKey}).

\subsection{$I_2(T_n,m,r)$}
By means of
\[
\int_{-1}^{1}\hypsng\frac{T_n(s)(1-s^2)^{m-\frac{1}{2}}}{(s-r)^2}ds =
\frac{d}{dr}\int_{-1}^{1}\cpv\frac{T_n(s)(1-s^2)^{m-\frac{1}{2}}}{s-r}ds
\ ,
\]
one readily obtains:
\begin{itemize}
\item{$I_2(T_n,0,r)$}~\cite{KayaErdo87}:
\begin{eqnarray}
\lefteqn{\int_{-1}^{1}\hypsng\frac{T_n(s)}{\sqrt{1-s^2}(s-r)^2}ds
  = \pi\frac{dU_{n-1}(r)}{dr}} \nonumber \\
 & \stackrel{(\ref{derv2Cheby})}{=} &
\frac{\pi}{1-r^2}\left[\frac{n+1}{2}U_{n-2}(r) -
\frac{n-1}{2}U_{n}(r)\right] \ ,\ \ n\geq 2
\end{eqnarray}

\item{$I_2(T_n,1,r)$:}
\begin{eqnarray}
\lefteqn{\int_{-1}^{1}\hypsng\frac{T_n(s)\sqrt{1-s^2}}{(s-r)^2}ds
  = \frac{\pi}{2}\frac{d}{dr}[T_{n-1}(r) - T_{n+1}(r)] } \nonumber \\
 & \stackrel{(\ref{derv1Cheby})}{=} &
\frac{\pi}{2}\left[(n-1)U_{n-2}(r) - (n+1)U_{n}(r)\right] \ ,\ \ n\geq 2
\end{eqnarray}

\item{$I_2(T_n,2,r)$:}
\begin{eqnarray}
\lefteqn{\int_{-1}^{1}\hypsng\frac{T_n(s)(1-s^2)^\frac{3}{2}}{(s-r)^2}ds
  = -\frac{\pi}{8}\frac{d}{dr}[T_{n-3}(r) - 3T_{n-1}(r) + 3T_{n+1}(r)
  - T_{n+3}(r)] } \nonumber \\
 & \stackrel{(\ref{derv1Cheby})}{=} &
-\frac{\pi}{8}\left[\ (n-3)U_{n-4}(r) - 3(n-1)U_{n-2}(r) + 3(n+1)U_{n}(r)\ 
  - \ \ \ \ \  
  \right. \nonumber \\
 & \  & \ \ \left. \ \ \ \ \  (n+3)U_{n+2}(r)\ \right] \ ,\ \ \ n\geq 4
\end{eqnarray}

\item{$I_2(T_n,3,r)$:}
\begin{eqnarray}
\lefteqn{\int_{-1}^{1}\hypsng\frac{T_n(s)(1-s^2)^\frac{5}{2}}{(s-r)^2}ds
 }  \nonumber \\
 & = & \frac{\pi}{32}\frac{d}{dr}[T_{n-5}(r) - 5T_{n-3}(r) +
 10T_{n-1}(r) - 10T_{n+1}(r) + 5T_{n+3}(r) - T_{n+5}(r)]  \nonumber \\
 & \stackrel{(\ref{derv1Cheby})}{=} &
\frac{\pi}{32}\left[(n-5)U_{n-6}(r) - 5(n-3)U_{n-4}(r) + 10(n-1)U_{n-2}(r) -
  10(n+1)U_{n}(r)  \right. \nonumber \\
 & \  & \ \ \left. \ \ \ \ \ + 5(n+3)U_{n+2}(r) - (n+5)U_{n+4}(r)\
 \right] \ ,\ \ \ n\geq 6
\end{eqnarray}

\item{$I_2(T_n,m,r)$, where $m\ge 1$, and $n\ge 2m+1$:}
\begin{equation}
\boxed{
\begin{array}{l}\mbox{
{\Large{
$\int_{-1}^{1}\hypsng\frac{T_n(s)(1-s^2)^{m-1/2}}{(s-r)^2}$}}}
ds =
\\ 
\pi(-1)^{m+1}\left(\frac{1}{2}\right)^{2m-1}\sum_{j=0}^{2m-1}(-1)^{j}\left(
\hspace{-.05in}\begin{array}{c}^{2m-1}
\\ _{j} \end{array} \hspace{-.05in}\right)(n+1-2m+2j)U_{n-2m+2j}(r)
\end{array}
}\label{boxIII}
\end{equation}


\end{itemize}

\subsection{$I_2(U_n,m,r)$}
The following equality
\[
\int_{-1}^{1}\hypsng\frac{U_n(s)(1-s^2)^{m-\frac{1}{2}}}{(s-r)^2}ds =
\frac{d}{dr}\int_{-1}^{1}\cpv\frac{U_n(s)(1-s^2)^{m-\frac{1}{2}}}{s-r}ds
\]
leads to:
\begin{itemize}
\item{$I_2(U_n,1,r)$}~\cite{KayaErdo87}:
\begin{equation}
\int_{-1}^{1}\hypsng\frac{U_n(s)\sqrt{1-s^2}}{(s-r)^2}ds
  = -\pi\frac{dT_{n+1}(r)}{dr}
 \stackrel{(\ref{derv1Cheby})}{=} -\pi(n+1)U_{n}(r) \ ,\ \ n\geq 0\ ,
\end{equation}
which is the same as (\ref{int2UCheby}).

\item{$I_2(U_n,2,r)$:}
\begin{eqnarray}
\lefteqn{\int_{-1}^{1}\hypsng\frac{U_n(s)(1-s^2)^\frac{3}{2}}{(s-r)^2}ds
  = \frac{\pi}{4}\frac{d}{dr}[T_{n-1}(r) - 2T_{n+1}(r) + T_{n+3}(r)] } \nonumber \\
 & \stackrel{(\ref{derv1Cheby})}{=} &
\frac{\pi}{4}\left[\ (n-1)U_{n-2}(r) - 2(n+1)U_{n}(r) + (n+3)U_{n+2}(r)
  \right] \ ,\ n\geq 2
\end{eqnarray}

\item{$I_2(U_n,3,r)$:}
\begin{eqnarray}
\lefteqn{\int_{-1}^{1}\hypsng\frac{U_n(s)(1-s^2)^\frac{5}{2}}{(s-r)^2}ds
 }  \nonumber \\
 & = & -\frac{\pi}{16}\frac{d}{dr}[T_{n-3}(r) - 4T_{n-1}(r) + 6T_{n+1}(r)
  - 4T_{n+3}(r) + T_{n+5}(r)]  \nonumber \\
 & \stackrel{(\ref{derv1Cheby})}{=} &
-\frac{\pi}{16}\left[(n-3)U_{n-4}(r) - 4(n-1)U_{n-2}(r) +
  6(n+1)U_{n}(r) - 4(n+3)U_{n+2}(r) \right. \nonumber \\
 & \  & \ \ \left. \ \ \ \ \  + (n+5)U_{n+4}(r)\ \right] \ ,\ \ \ n\geq 4
\end{eqnarray}

\item{$I_2(U_n,m,r)$, where $m\ge 2$, and $n\ge 2m-1$:}
\begin{equation}
\boxed{
\begin{array}{l}\mbox{\Large
$\int_{-1}^{1}\hypsng\frac{U_n(s)(1-s^2)^{m-\frac{1}{2}}}{(s-r)^2}$}
ds=
\\
\pi(-1)^{m}\left(\frac{1}{2}\right)^{2m-2}\sum_{j=0}^{2m-2}(-1)^{j}\left(
\hspace{-.05in}\begin{array}{c} ^{2m-2}
\\ _{j} \end{array} \hspace{-.05in}\right)(n+3-2m+2j)U_{n+2-2m+2j}(r)
\end{array}
}\label{boxIV}
\end{equation}

\end{itemize}

\subsection{$I_3(T_n,m,r)$}
By means of
\[
\int_{-1}^{1}\hypsng\frac{T_n(s)(1-s^2)^{m-\frac{1}{2}}}{(s-r)^3}ds =
\frac{1}{2}\frac{d}{dr}\int_{-1}^{1}\hypsng\frac{T_n(s)(1-s^2)^{m-\frac{1}{2}}}
{(s-r)^2}ds\ ,
\]
one obtains:

\begin{itemize}

\item{$I_3(T_n,0,r)$:}
\begin{eqnarray}
\lefteqn{\int_{-1}^{1} \hypsng
\frac{T_n(s)\sqrt{1-s^2}}{(s-r)^3}ds =
\frac{\pi}{8(1-r^2)^2}\mbox{\Large{[}}(n+1)(n+2)U_{n-3}(r) } \nonumber \\
 & \  & -2(n^2-3)U_{n-1}(r) + (n-1)^2U_{n+1}(r)\mbox{\Large{]}}\ , \ \ \ n \geq 3
\end{eqnarray}

\item{$I_3(T_n,1,r)$:}
\begin{eqnarray}
\lefteqn{\int_{-1}^{1} \hypsng
\frac{T_n(s)\sqrt{1-s^2}}{(s-r)^3}ds =
\frac{\pi}{8(1-r^2)}\mbox{\Large{[}}(n^2-n)U_{n-3}(r) } \nonumber \\
 & - & (2n^2+2)U_{n-1}(r) + (n^2+n)U_{n+1}(r)\mbox{\Large{]}} \ , \ \ \ n \geq 3
\end{eqnarray}

\item{$I_3(T_n,2,r)$:}
\begin{eqnarray} 
\lefteqn{\int_{-1}^{1} \hypsng
 \frac{T_n(s)(1-s^2)^{\frac{3}{2}}}{(s-r)^3}ds=\frac{\pi}{32(1-r^2)}\mbox{\Large{\{}}-(n+3)(n+2)U_{n+3}(r) }\nonumber \\
 & + & [(n+3)(n+4)+3n(n+1)]U_{n+1}(r)-[3(n+1)(n+2)+3(n-1)(n-2)]U_{n-1}(r) \nonumber \\
 & + & [3n(n-1)+(n-3)(n-4)]U_{n-3}(r) - (n-3)(n-2)U_{n-5}(r)\mbox{\Large{\}}}\ , \ \ n\geq 5
\label{hyperCube}
\end{eqnarray}

\item{$I_3(T_n,3,r)$:}
\begin{eqnarray} 
\lefteqn{\int_{-1}^{1} \hypsng
 \frac{T_n(s)(1-s^2)^{\frac{5}{2}}}{(s-r)^3}ds=\frac{\pi}{128(1-r^2)}\mbox{\Large{[}}(n^2+9n+20)U_{n+5}(r)  } \nonumber \\
 & - &
 6(n^2+6n+10)U_{n+3}(r)+15(n^2+3n+4)U_{n+1}(r)  \nonumber \\
 & - & 20(n^2+2)U_{n-1}(r) + 15(n^2-3n+4)U_{n-3}(r)-
 6(n^2-6n+10)U_{n-5}(r) \nonumber \\
 & + & (n^2-9n+2)U_{n-7}(r)\mbox{\Large{]}}\ ,\ \ \ n\geq 7
\label{myRefExp3}
\end{eqnarray}

\item{$I_3(T_n,m,r)$, where $m\ge 1$, and $n\ge 2m+2$:}
\begin{equation}
\boxed{\begin{array}{l}\mbox{\Large
$\int_{-1}^{1}\hypsng\frac{T_n(s)(1-s^2)^{m-\frac{1}{2}}}{(s-r)^3}$}ds=
 \\ 
(-1)^{m+1}\left(\frac{1}{2}\right)^{2m+1}\frac{\pi}{1-r^2}\sum_{j=0}^{2m-1}(-1)^{j}\left(
\hspace{-.05in}\begin{array}{c} ^{2m-1}
\\ _{j} \end{array} \hspace{-.05in}\right)(n+1-2m+2j)\times \\ 
\left[(n+2-2m+2j)U_{n-1-2m+2j}(r)-(n-2m+2j)U_{n+1-2m+2j}(r)
\right]
\end{array}
}\label{boxV}
\end{equation}


\end{itemize}

\subsection{$I_3(U_n,m,r)$}
By means of
\[
\int_{-1}^{1}\hypsng\frac{U_n(s)(1-s^2)^{m-\frac{1}{2}}}{(s-r)^3}ds =
\frac{1}{2}\frac{d}{dr}\int_{-1}^{1}\hypsng\frac{U_n(s)(1-s^2)^{m-\frac{1}{2}}}
{(s-r)^2}ds\ ,
\]
one gets:
\begin{itemize}
\item{$I_3(U_n,1,r)$:}
\begin{eqnarray}
\lefteqn{\int_{-1}^{1} \hypsng
\frac{U_n(s)\sqrt{1-s^2}}{(s-r)^3}ds } \nonumber \\
& = &
\frac{\pi}{4(1-r^2)}\left\{
  -  (2n^2+3n+2)U_{n-1}(r) + (n^2+n)U_{n+1}(r)\right\}\ , \ n \geq 1
\end{eqnarray}

\item{$I_3(U_n,2,r)$:}
\begin{eqnarray} 
\lefteqn{\int_{-1}^{1} \hypsng
 \frac{U_n(s)(1-s^2)^{\frac{3}{2}}}{(s-r)^3}ds=\frac{\pi}{16(1-r^2)}\mbox{\Large{[}}-(n^2+5n+6)U_{n+3}(r)
 } \nonumber \\
 & + & (3n^2+9n+12)U_{n+1}(r) - (3n^2+3n+6)U_{n-1}(r) \nonumber \\
 & + & (n^2-n)U_{n-3}(r)\mbox{\Large{]}}\ , \ \ \ n\geq 3
\end{eqnarray}


\item{$I_3(U_n,3,r)$:}
\begin{eqnarray} 
\lefteqn{\int_{-1}^{1} \hypsng
 \frac{U_n(s)(1-s^2)^{\frac{5}{2}}}{(s-r)^3}ds=\frac{\pi}{64(1-r^2)}\mbox{\Large{[}}(n^2+9n+20)U_{n+5}(r)
 } \nonumber \\
 &-&(5n^2+31n+54)U_{n+3}(r)+(10n^2+34n+48)U_{n+1}(r)-(10n^2+6n+20)U_{n-1}(r)\nonumber \\
 &+&(5n^2-11n+12)U_{n-3}(r)-(n^2-5n+6)U_{n-5}(r)\mbox{\Large{]}}\ , \ \ \ n\geq 5
\end{eqnarray}

\item{$I_3(U_n,m,r)$, where $m\ge 2$, and $n\ge 2m$:}
\begin{equation}
\boxed{
\begin{array}{l}\mbox{\Large
$\int_{-1}^{1}\hypsng\frac{U_n(s)(1-s^2)^{m-\frac{1}{2}}}{(s-r)^3}$}
ds= \\ 
(-1)^{m}\left(\frac{1}{2}\right)^{2m}\frac{\pi}{1-r^2}\sum_{j=0}^{2m-2}(-1)^{j}\left(
\hspace{-.05in}\begin{array}{c} ^{2m-2}
\\ _{j} \end{array} \hspace{-.05in}\right)(n+3-2m+2j)\times \\ 
\left[(n+4-2m+2j)U_{n+1-2m+2j}(r)-(n+2-2m+2j)U_{n+3-2m+2j}(r)
\right]
\end{array}
}\label{boxVI}
\end{equation}


\end{itemize}

\subsection{$I_4(T_n,m,r)$}
By means of
\[
\int_{-1}^{1}\hypsng\frac{T_n(s)(1-s^2)^{m-\frac{1}{2}}}{(s-r)^4}ds =
\frac{1}{3}\frac{d}{dr}\int_{-1}^{1}\hypsng\frac{T_n(s)(1-s^2)^{m-\frac{1}{2}}}
{(s-r)^3}ds\ ,
\]
one reaches the following results:
\begin{itemize}
\item{$I_4(T_n,0,r)$:}
\begin{eqnarray}
\lefteqn{\int_{-1}^{1} \hypsng
 \frac{T_n(s)}{(s-r)^4\sqrt{1-s^2}}ds=\frac{\pi}{48(1-r^2)^3}\mbox{\Large{[}}
 (n^3+6n^2+11n+6)U_{n-4}(r)  } \nonumber \\
 &\  &  -(3n^3+6n^2-25n-44)U_{n-2}(r) +(3n^3-5n^2-19n+37)U_{n}(r) \nonumber \\
 &\  & -(n^3-5n^2+7n-3)U_{n+2}(r) \mbox{\Large{]}}\  ,\ \ \ \ \ \ \ \  n\geq 4
\end{eqnarray}

\item{$I_4(T_n,1,r)$:}
\begin{eqnarray}
\lefteqn{\int_{-1}^{1} \hypsng
 \frac{T_n(s)\sqrt{1-s^2}}{(s-r)^4}ds=\frac{\pi}{48(1-r^2)^2}\mbox{\Large{[}}
  } \nonumber \\
 &\  &  (n^3-n)U_{n-4}(r)-(3n^3+9n+12)U_{n-2}(r) +(3n^3+9n-12)U_{n}(r) \nonumber \\
 &\  & -(n^3-n)U_{n+2}(r) \mbox{\Large{]}}\  ,\ \ \ \ \ \ \ \  n\geq 4
\end{eqnarray}

\item{$I_4(T_n,2,r)$:}
\begin{eqnarray} 
\lefteqn{\int_{-1}^{1} \hypsng
 \frac{T_n(s)(1-s^2)^{\frac{3}{2}}}{(s-r)^4}ds=\frac{\pi}{192(1-r^2)^2}\mbox{\Large{[}}
 (n^3+6n^2+11n+6)U_{n+4}(r) } \nonumber \\
 & - &
 (5n^3+18n^2+43n+30)U_{n+2}(r) + (10n^3+12n^2+134n-36)U_{n}(r)\nonumber \\
 & - &
 (10n^3-12n^2+134n+36)U_{n-2}(r)+ (5n^3-18n^2+43n-30)U_{n-4}(r) \nonumber \\
 & - &
 (n^3-6n^2+11n-6)U_{n-6}(r)\mbox{\Large{]}},\ \ \ \ \ \ n\geq 6
\end{eqnarray}

\item{$I_4(T_n,3,r)$:}
\begin{eqnarray} 
\lefteqn{\int_{-1}^{1} \hypsng
 \frac{T_n(s)(1-s^2)^{\frac{5}{2}}}{(s-r)^4}ds=\frac{\pi}{384(1-r^2)^2}\mbox{\Large{[}}
 -(\frac{1}{2}n^3+6n^2+\frac{47}{2}n+30)U_{n+6}(r) } \nonumber \\
 & + &
 (\frac{7}{2}n^3+30n^2+\frac{197}{2}n+120)U_{n+4}(r) -
 (\frac{21}{2}n^3+54n^2+\frac{327}{2}n+180)U_{n+2}(r) \nonumber \\
 & + &
 (\frac{35}{2}n^3+30n^2+\frac{325}{2}n+90)U_{n}(r) -
 (\frac{35}{2}n^3-30n^2+\frac{325}{2}n-90)U_{n-2}(r) \nonumber \\
 & + & (\frac{21}{2}n^3-54n^2+\frac{327}{2}n-180)U_{n-4}(r) -
 (\frac{7}{2}n^3-30n^2+\frac{197}{2}n-120)U_{n-6}(r) \nonumber \\
 & + &
 (\frac{1}{2}n^3-6n^2+\frac{47}{2}n-30)U_{n-8}(r)\mbox{\Large{]}},\ \ \ \   n\geq 8
\end{eqnarray}

\item{$I_4(T_n,m,r)$, where $m\ge 1$, and $n\ge 2m+3$:}
\begin{equation}
\boxed{
\begin{array}{l}\mbox{\Large
$\int_{-1}^{1}\hypsng\frac{T_n(s)(1-s^2)^{m-\frac{1}{2}}}{(s-r)^4}$}
ds = \\ 
(-1)^{m+1}\left(\frac{1}{2}\right)^{2m+2}\frac{1}{3}\frac{\pi}{(1-r^2)^2}
\sum_{j=0}^{2m-1}(-1)^{j}\left(
\hspace{-.05in}\begin{array}{c} ^{2m-1}
\\ _{j} \end{array} \hspace{-.05in}\right)(n+1-2m+2j)\times  \\ 
 \mbox{\Large{\{}} 
[(n+2-2m+2j)(n+3-2m+2j)]U_{n-2-2m+2j}(r) - \\
 \ \ [2(n-2m+2j)^2+4(n-2m+2j)-6]U_{n-2m+2j}(r) + \\
 \ \ [(n-2m+2j)(n-1-2m+2j)]U_{n+2-2m+2j}(r)\mbox{\Large{\}}}
\end{array}
}\label{boxVII}
\end{equation}


\end{itemize}

\subsection{$I_4(U_n,m,r)$:}
By means of
\[
\int_{-1}^{1}\hypsng\frac{U_n(s)(1-s^2)^{m-\frac{1}{2}}}{(s-r)^4}ds =
\frac{1}{3}\frac{d}{dr}\int_{-1}^{1}\hypsng\frac{U_n(s)(1-s^2)^{m-\frac{1}{2}}}
{(s-r)^3}ds\ ,
\]
one obtains
\begin{itemize}
\item{$I_4(U_n,1,r)$:}
\begin{eqnarray}
\lefteqn{\int_{-1}^{1} \hypsng
 \frac{U_n(s)(1-s^2)^{\frac{1}{2}}}{(s-r)^4}ds=\frac{\pi}{24(1-r^2)^2}\mbox{\Large{[}}
  } \nonumber \\
 &\  &  -(2n^3+9n^2+11n+6)U_{n-2}(r) +(3n^3+3n^2-2n-6)U_{n}(r) \nonumber \\
 &\  & -(n^3-n)U_{n+2}(r) \mbox{\Large{]}}\  ,\ \ \ \ \ \ \ \  n\geq 2
\end{eqnarray}

\item{$I_4(U_n,2,r)$:}
\begin{eqnarray}
\lefteqn{\int_{-1}^{1} \hypsng
 \frac{U_n(s)(1-s^2)^{\frac{3}{2}}}{(s-r)^4}ds=\frac{\pi}{96(1-r^2)^2}\mbox{\Large{[}}
  (n^3+6n^2+11n+6)U_{n+4}(r)} \nonumber \\
 &\  &  -(4n^3+18n^2+44n+30)U_{n+2}(r) +(6n^3+18n^2+54n+42)U_{n}(r) \nonumber \\
 &\  & -(4n^3+6n^2+20n+18)U_{n-2}(r)+ (n^3-n)U_{n-4}(r)\mbox{\Large{]}}\  ,\ \ \ \ \ \  n\geq 4
\end{eqnarray}

\item{$I_4(U_n,3,r)$:}
\begin{eqnarray} 
\lefteqn{\int_{-1}^{1} \hypsng
 \frac{U_n(s)(1-s^2)^{\frac{5}{2}}}{(s-r)^4}ds=\frac{\pi}{192(1-r^2)^2}\mbox{\Large{[}}
 -(\frac{1}{2}n^3+6n^2+\frac{47}{2}n+320)U_{n+6}(r) } \nonumber \\
 & + &
 (3n^3+27n^2+93n+117)U_{n+4}(r) - (\frac{15}{2}n^3+45n^2+\frac{285}{2}n+165)U_{n+2}(r)
 \nonumber \\
 & + &
 (10n^3+30n^2+110n+90)U_{n}(r) - (\frac{15}{2}n^3+\frac{105}{2}n)U_{n-2}(r)
 \nonumber \\
 & + & (3n^3-9n^2+21n-15)U_{n-4}(r) \nonumber \\
 & - & (\frac{1}{2}n^3+3n^2+\frac{11}{2}n-3)U_{n-6}(r)\mbox{\Large{]}},\ \ \ \ n\geq 6
\end{eqnarray}

\item{$I_4(U_n,m,r)$, where $m\ge 2$, and $n\ge 2m+1$:}
\begin{equation}
\boxed{
\begin{array}{l}\mbox{\Large
$\int_{-1}^{1}\hypsng\frac{U_n(s)(1-s^2)^{m-\frac{1}{2}}}{(s-r)^4}$}
ds = \\ 
(-1)^{m}\left(\frac{1}{2}\right)^{2m+1}\frac{1}{3}\frac{\pi}{(1-r^2)^2}
\sum_{j=0}^{2m-2}(-1)^{j}\left(
\hspace{-.05in}\begin{array}{c} ^{2m-2}
\\ _{j} \end{array} \hspace{-.05in}\right)(n+3-2m+2j)\times  \\ 
 \mbox{\Large{\{}} 
[(n+4-2m+2j)(n+5-2m+2j)]U_{n-2m+2j}(r) - \\
 \ \ [2(n-2m+2j)^2 + 10(n-2m+2j) +10]U_{n+2-2m+2j}(r) + \\
 \ \ [(n+2-2m+2j)(n-1-2m+2j)]U_{n+4-2m+2j}(r)\mbox{\Large{\}}}
\end{array}
}\label{boxVIII}
\end{equation}


\end{itemize}

\section{Evaluation of Stress Intensity Factors (SIFs)}
An important task is to evaluate the stress intensity factors
(SIFs) at both crack tips, since the propagation of a crack starts
around its tips.  In mode III crack problems, SIFs can be calculated from
\begin{eqnarray}
K_{III}(d)
 & = & \lim_{x\rightarrow d^+}\sqrt{2\pi(x-d)}\sigma_{yz}(x,0) \ \ , \ \
 \ ( x > d ) \label{sifDefP}
\end{eqnarray}
and
\begin{eqnarray}
K_{III}(c)
 & = & \lim_{x\rightarrow c^-}\sqrt{2\pi(c-x)}\sigma_{yz}(x,0) \ \ , \ \
 \ ( x < c ) \ .\label{sifDefN}
\end{eqnarray}

Note that the limit is taken from outside of the crack surfaces
and towards
both tips. Usually the left hand side of integral equation
(\ref{int_xt2}) is the
expression for $\sigma_{yz}(x,0)$ which is valid for $x$ is inside the
crack surfaces $(c, d)$ as well as outside of $(c, d)$.  
Thus to calculate SIFs, the key is to evaluate the
following integrals which are obtained after proper normalization and
the change of variables described in equation (\ref{chngVar}), 
\begin{equation}
S_{\alpha}(T_n,m,r) = 
\int_{-1}^{1}\frac{T_n(s)(1-s^2)^{m-(1/2)}}{(s-r)^{\alpha}}ds
\ , \ \ \ r \notin (-1, 1) \label{sifTintNorm}
\end{equation}
and
\begin{equation}
S_{\alpha}(U_n,m,r) = 
\int_{-1}^{1}\frac{U_n(s)(1-s^2)^{m-(1/2)}}{(s-r)^{\alpha}}ds
\ , \ \ \ r \notin (-1, 1)  \ . \label{sifUintNorm}
\end{equation}
Note that the above integrals are not singular as
$x\neq t$ for $t \in (c, d)$ and $x
\notin (c, d)$.


\

The strategy to evaluate $S_{\alpha}(T_n,m,r)$ and
$S_{\alpha}(U_n,m,r)$ for general integers $\alpha$ (positive)
and $m$ (non-negative) is similar to the process for evaluating
$I_{\alpha}(T_n,m,r)$ and $I_{\alpha}(U_n,m,r)$.  It consists of 
evaluating the integrals $S_{1}(T_n,m,r)$ and $S_{1}(U_n,m,r)$ by
means of the reduction procedure described in {\textsf{Section 4.3}},
and taking differentiation (with respect to $r$) to obtain
$S_{\alpha}(T_n,m,r)$ and $S_{\alpha}(U_n,m,r)$ for $\alpha \geq 2$.
The relevant derivations are provided below where  the range of $r$ is
restricted to $|r|>1$. These formulas are used in calculating SIFs
for the examples presented in \textsf{Section 7}

\subsection{$S_{1}(T_n,m,r)$ and $S_{1}(U_n,m,r)$}

\begin{itemize}

\item $S_{1}(T_n,0,r)$ (This is a well known integral~\cite{ErdoOztu92}):

\begin{equation}
\int_{-1}^{1}\frac{T_n(s)}{(s-r)\sqrt{1-s^2}} ds =
-\pi\frac{\left(r-\sqrt{r^2-1}|r|/r\right)^n}{\sqrt{r^2-1}|r|/r} \
, \quad n \geq 0\ .  \label{out1Tn0}
\end{equation}

\item $S_{1}(T_n,1,r)$:
\begin{eqnarray}
\lefteqn{\int_{-1}^{1}\frac{T_n(s)\sqrt{1-s^2}}{s-r}\ ds} \nonumber \\
 & \stackrel{(\ref{recur3Cheby})}{=}&
\frac{1}{2}\int_{-1}^{1}\frac{U_{n}(s)\sqrt{1-s^2}}{s-r} ds -
\frac{1}{2}\int_{-1}^{1}\frac{U_{n-2}(s)\sqrt{1-s^2}}{s-r}\ ds \nonumber \\
 & \stackrel{(\ref{out1Un1})}{=}& \pi\frac{|r|}{r}\sqrt{r^2-1}\left(r -
 \frac{|r|}{r}\sqrt{r^2-1}\right)^{n} , \quad n \geq 2\ .
\label{out1Tn1}
\end{eqnarray}

\item $S_{1}(T_0,2,r)$:
\begin{equation}
\int_{-1}^{1}\frac{T_0(s)(1-s^2)^{\frac{3}{2}}}{s-r} ds =
\pi(r^2-1)\left(r - \frac{|r|}{r}\sqrt{r^2-1}\right)\ .
\label{cauchy0}
\end{equation}

\item $S_{1}(T_1,2,r)$:
\begin{equation}
\int_{-1}^{1}\frac{T_1(s)(1-s^2)^{\frac{3}{2}}}{s-r}\
ds =
\frac{\pi}{2}(r^2-1)\left(r - \frac{|r|}{r}\sqrt{r^2-1}\right)^2\ .
\label{cauchy1}
\end{equation}

\item $S_{1}(T_n,2,r)$:
\begin{equation}
\int_{-1}^{1}\frac{T_n(s)(1-s^2)^{\frac{3}{2}}}{s-r} ds
\stackrel{(\ref{recur3Cheby}),(\ref{out1Un2})}{\mbox{\Large =}}
-\frac{\pi|r|}{r}(r^2-1)^{\frac{3}{2}}\left(r - \frac{|r|}{r}\sqrt{r^2-1}\right)
^n \
, \ n \geq 2 \ . \label{cauchyGnrL}
\end{equation}

\item $S_{1}(U_n,1,r)$:
\begin{eqnarray}
\lefteqn{\int_{-1}^{1}\frac{U_n(s)\sqrt{1-s^2}}{s-r}\ ds} \nonumber \\
 & \stackrel{(\ref{recur5Cheby})}{=}&
\frac{1}{2}\int_{-1}^{1}\frac{T_{n}(s)}{(s-r)\sqrt{1-s^2}} ds -
\frac{1}{2}\int_{-1}^{1}\frac{T_{n+2}(s)}{(s-r)\sqrt{1-s^2}}\ ds \nonumber \\
 & \stackrel{(\ref{out1Tn0})}{=}& -\pi\left(r -
 \frac{|r|}{r}\sqrt{r^2-1}\right)^{n+1} , \quad n \geq 0\ .
\label{out1Un1}
\end{eqnarray}

\item $S_{1}(U_n,2,r)$:
\begin{eqnarray}
\lefteqn{\int_{-1}^{1}\frac{U_n(s)(1-s^2)^{\frac{3}{2}}}{s-r}\ ds} \nonumber \\
 & \stackrel{(\ref{recur5Cheby})}{=}&
\frac{1}{2}\int_{-1}^{1}\frac{T_{n}(s)\sqrt{1-s^2}}{s-r} ds -
\frac{1}{2}\int_{-1}^{1}\frac{T_{n+2}(s)\sqrt{1-s^2}}{s-r}\ ds \nonumber \\
 & \stackrel{(\ref{out1Tn1})}{=}& \pi(r^2-1)\left(r -
 \frac{|r|}{r}\sqrt{r^2-1}\right)^{n+1} , \quad n \geq 2\ .
\label{out1Un2}
\end{eqnarray}

\end{itemize}


\

The formulas for $S_{1}(T_n,m,r)$ and
$S_{1}(U_n,m,r)$ with general $m$ can be deduced by
the same procedure described in \textsf{Figure~\ref{fig0}},
and are listed below.
\begin{eqnarray}
\lefteqn{S_{1}(T_n,m,r)  =
\int_{-1}^{1}\frac{T_n(s)(1-s^2)^{m-1/2}}{s-r} ds}
\nonumber \\ & = &
\pi(-1)^{m+1}\frac{|r|}{r}(r^2-1)^{m-1/2}
\left(r - \frac{|r|}{r}\sqrt{r^2-1}\right)^n \
, \ m \geq 0 \mbox{ and }\ n \geq 2m\ .
\label{mCauchyTn}
\end{eqnarray}

\ 

\begin{eqnarray}
\lefteqn{S_{1}(U_n,m,r)  =
\int_{-1}^{1}\frac{U_n(s)(1-s^2)^{m-1/2}}{s-r} ds}
\nonumber \\ & = &
\pi(-1)^{m}(r^2-1)^{m-1}
\left(r - \frac{|r|}{r}\sqrt{r^2-1}\right)^n \
, \ m \geq 1 \mbox{ and }\ n \geq 2m-2\ .
\label{mCauchyUn}
\end{eqnarray}

\subsection{$S_{2}(T_n,m,r)$ and $S_{2}(U_n,m,r)$}

Differentiating (with respect to $r$) the formulas for
$S_{1}(T_n,m,r)$ and $S_{1}(U_n,m,r)$, 
we obtain the formulas for $S_{2}(T_n,m,r)$ and $S_{2}(U_n,m,r)$.
\begin{eqnarray}
\lefteqn{\int_{-1}^{1}\frac{U_n(s)(1-s^2)^{\frac{3}{2}}}{(s-r)^2}\
ds } \nonumber \\
 & = &
-\pi(n+1)\frac{|r|}{r}\sqrt{r^2-1}\left(r - \frac{|r|}{r}\sqrt{r^2-1}\right)^{n+1} +
2\pi r\left(r - \frac{|r|}{r}\sqrt{r^2-1}\right)^{n+1}
, \ n \geq 0
\label{outsideInt2}
\end{eqnarray}
and
\begin{equation}
\int_{-1}^{1}\frac{T_n(s)(1-s^2)^{\frac{3}{2}}}{s-r} ds =
-\frac{\pi|r|}{r}(r^2-1)^{\frac{3}{2}}\left(r - \frac{|r|}{r}\sqrt{r^2-1}\right)^n \
, \ n \geq 2  \ .\label{cauchyn}
\end{equation}

\subsection{$S_{3}(T_n,m,r)$ and $S_{3}(U_n,m,r)$}
The following formulas are obtained by
differentiating twice (with respect to $r$) the corresponding
formulas obtained in {\sf Subsection 6.1}.  
\begin{eqnarray}
\lefteqn{\int_{-1}^{1}\frac{U_n(s)(1-s^2)^{\frac{3}{2}}}{(s-r)^3} ds
=}
 \nonumber \\ & \ &
\frac{\pi}{2}\left[(n^2+2n+3)-3(n+1)\frac{|r|}{\sqrt{r^2-1}}\right]
\left(r - \frac{|r|}{r}\sqrt{r^2-1}\right)^{n+1}, \ n \geq 0\ .
\label{outsideInt3}
\end{eqnarray}

\begin{equation}
\int_{-1}^{1}\frac{T_1(s)(1-s^2)^{\frac{3}{2}}}{(s-r)^3}ds
=\frac{3\pi}{2}\left(r - \frac{|r|}{r}\sqrt{r^2-1}\right)^2
\left(1 - \frac{|r|}{\sqrt{r^2-1}}\right)\ .
\end{equation}

\begin{equation}
\int_{-1}^{1}\frac{T_0(s)(1-s^2)^{\frac{3}{2}}}{(s-r)^3}ds
=\frac{3\pi}{2}\left(r - \frac{|r|}{r}\sqrt{r^2-1}\right)
\left(1 - \frac{|r|}{\sqrt{r^2-1}}\right)\ .
\end{equation}

\begin{eqnarray}
\lefteqn{\int_{-1}^{1}\frac{T_n(s)(1-s^2)^{\frac{3}{2}}}{(s-r)^3}\
ds } \nonumber \\
 & = & \frac{\pi}{4}\left\{\left(r - \frac{|r|}{r}\sqrt{r^2-1}\right)^{n+1}\left[
(n^2+2n+3)-3(n+1)\frac{|r|}{\sqrt{r^2-1}}\right]\right. \nonumber \\
 & \ & -\left.\left(r - \frac{|r|}{r}\sqrt{r^2-1}\right)^{n-1}\left[
(n^2-2n+3)-3(n-1)\frac{|r|}{\sqrt{r^2-1}}\right]\right\} \ , \ n \geq 2\ .
\end{eqnarray}


\section{\large{Examples}}
\label{exam}
Three examples are presented here, which emphasize various aspects of
Fredholm singular integral equation formulations and their linkage to
fracture mechanics.  These examples are:

1. Internal mode I crack in an infinite strip.

2. Mode III crack problem in nonhomogeneous materials.

3. Gradient elasticity theory applied to a mode III crack.

\

The first and last examples consider homogeneous materials, and the second
example considers nonhomogeneous materials, which has relevant
applications to the field of functionally graded
materials~\cite{ChanPaulFann99a,PaulFannChanFGM98}.
The first two examples are
from classical elasticity, and the last one is from
gradient elasticity theory.  The first example involves
mode I cracks and the last two examples involve mode III cracks.  All the
examples are formulated by using hypersingular integral equations.  For the
first two examples the order of singularity $\alpha$ is 2, and for the last
example $\alpha$ is 3.  A detailed comparison between $U_n$ and $T_n$ representations
is given in the first example.  A discussion on the
influence of the density function on the order of singularity of the
integral equation is
presented in the second and third examples.  The description of the examples is
summarized in {\sf Table \ref{tab1}}

\begin{table}[htbp]
\centering
\caption{Description of the examples.}
\begin{tabular}{||l||c|c|c||}   \hline \hline
{\bf Description} & {\bf Example 1} & {\bf Example 2} & {\bf Example 3} \\ \hline \hline
Homogeneous material & $\surd$ & \  & $\surd$ \\ \hline
Nonomogeneous material & \  & $\surd$ & \ \\ \hline
Classical elasticity & $\surd$ & $\surd$ & \  \\ \hline
Gradient elasticity & \  & \  & $\surd$ \\ \hline
Crack mode & I & III & III \\ \hline
Density function & displacement ($v$) & displacement ($w$) & slope ($\phi$) \\ \hline
Degree of singularity $\alpha$ & 2 & 2 & 3 \\ \hline
Weight function exponent, $m - (1/2)$ & 1/2 & 1/2 & 3/2 \\ \hline
Representation & $U_n$, $T_n$ & $U_n$, $T_n$ & $T_n$ \\ \hline
\end{tabular}
\label{tab1}
\end{table}

\subsection{\normalsize{Internal Mode I Crack in an Infinite
Strip~\cite{KayaErdo87}}}
\begin{figure}[ht]
\centerline{\epsfxsize= 12.72 cm  \epsffile{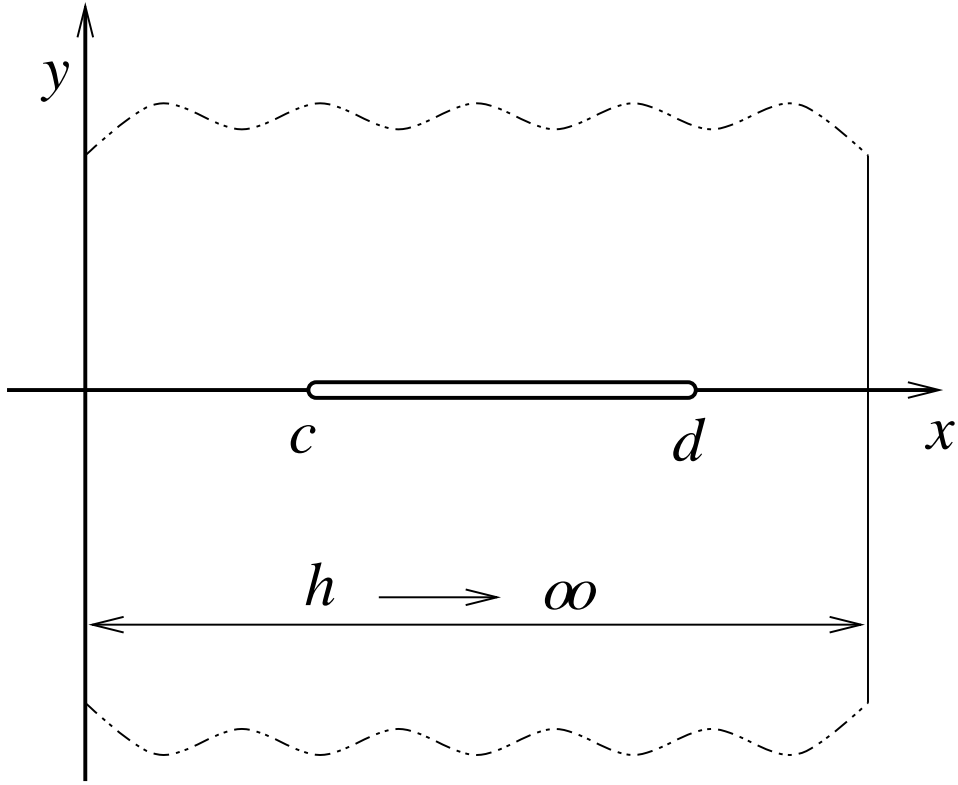}}
 \caption{A mode I crack in an infinite strip.}
 \label{fig1}
\end{figure}

Consider a crack in an infinite strip of homogeneous material, as illustrated by
\textsf{Figure~\ref{fig1}}. 
The governing partial differential equations ({\sf PDE}s) and boundary
conditions are:
\begin{equation}
\begin{array}{ll}
\nabla^{2}u(x,y) + \frac{2}{\kappa-1}\left(\frac{\partial^2
u(x,y)}{\partial x^2}+\frac{\partial^2u(x,y)}{\partial x\partial y} \right) = 0\ ,
& \quad -\infty<x, \ y<\infty\ , \\
\nabla^{2}v(x,y) + \frac{2}{\kappa-1}\left(\frac{\partial^2
v(x,y)}{\partial y^2}+\frac{\partial^2v(x,y)}{\partial x\partial y} \right) = 0\ ,
& \quad -\infty<x, \ y<\infty\ , \\
\sigma_{xx}(0, y) = \sigma_{xy}(0, y) = \sigma_{xx}(h, y) =
\sigma_{xy}(h, y) = 0\ , & \quad -\infty< y <\infty\ , \\
\sigma_{xy}(x, 0) = 0\ , & \quad 0< x < h\ , \\
\sigma_{yy}(x, 0) = -p(x) \ , & \quad x\in (c, d)\ , \\
v(x,0) = 0\ , & \quad x\notin [c, d]\ ,
\end{array}
\end{equation}
where $u$ and $v$ are the $x$ and $y$ components of the
displacement vector; $\sigma_{ij}$ is the stress
tensor; $\kappa$ is an elastic constant ($\kappa = 3-4\nu$ for plane
strain, $\kappa = (3-\nu)/(1+\nu)$ for plane stress, and $\nu$ is the
Poisson's ratio.)
This problem has been studied by Kaya and Erdogan~\cite{KayaErdo87} by
means of a $U_n$ representation, and it has also been used as a
benchmark
problem by Kabir {\it et. al.}~\cite{KabirETAL98}.  Here both $U_n$ and $T_n$ are
employed and compared.

\ 

The governing integral equation can be written in the form given by
equation (\ref{int_xt2}) as~\cite{KayaErdo87}
\begin{equation}
\int_{c}^{d}\hypsng\frac{\triangle v(t)}{(t-x)^2}dt +
\int_{c}^{d}k(x,t)\triangle v(t)dt =
-\pi\left(\frac{1+\kappa}{2\mu}\right)p(x)\ , \ \ c<x<d \ ,
\label{intKaya1}
\end{equation}
where the primary variable is the crack opening displacement
$\triangle v$ given by
\[ \triangle v(x) = v(x,0^+) - v(x,0^-) \ , \ \ c<x<d \ ,\]
and the kernel $k(x,t)$ is given in Kaya and
Erdogan~\cite{KayaErdo87}, equations (51) -- (54c), page 112.
It is worth noting that as $h \rightarrow \infty$ (see 
\textsf{Figure~\ref{fig1}}), the integral
equation for the half plane is recovered and the kernel $k(x,t)$ is
reduced to a much simpler form\footnote{To be consistent with the
notation adopted in this paper, we have used symbols different
from reference~\cite{KayaErdo87}.  For instance, in~\cite{KayaErdo87}
upper case $K(t,x)$ is
used for $k(t,x)$.}
\[ k(x,t) =
\frac{-1}{(t+x)^2}+\frac{12x}{(t+x)^3}-\frac{12x^2}{(t+x)^4} \ . \]
After normalization, the corresponding integral equation can be written
in a fashion similar to equation (\ref{intNorm}),
{\it i.e.}\footnote{Again, the notation is different from~\cite{KayaErdo87}.}
\begin{equation}
\int_{-1}^{1}\hypsng\frac{D(s)}{(s-r)^2}ds +
\int_{-1}^{1}\mathcal{K}(r,s)D(s)ds = P(r)\ , \ \ -1<r<1 \ ,
\label{intKaya2}
\end{equation}
where $D(s)$ is the unknown displacement function, the
regular kernel is
\[
\mathcal{K}(r,s) =
\frac{-1}{\left[(r+s)+2\left(\frac{d+c}{d-c}\right)\right]^2} +
\frac{12\left[s+\left(\frac{d+c}{d-c}\right)\right]}{\left[(r+s)+2\left(\frac{d+c}{d-c}\right)\right]^3}
-
\frac{12\left[s+\left(\frac{d+c}{d-c}\right)\right]^2}{\left[(r+s)+2\left(\frac{d+c}{d-c}\right)\right]^4}
\ ,
\]
and the loading function is
\[
P(r) =
-\pi\left(\frac{1+\kappa}{2\mu}\right)p\left(\left(\frac{d-c}{2}\right)s
+ \frac{d+c}{2} \right) \ .
\]

\begin{figure}[ht]
\centerline{\epsfxsize=14.72 cm  \epsffile{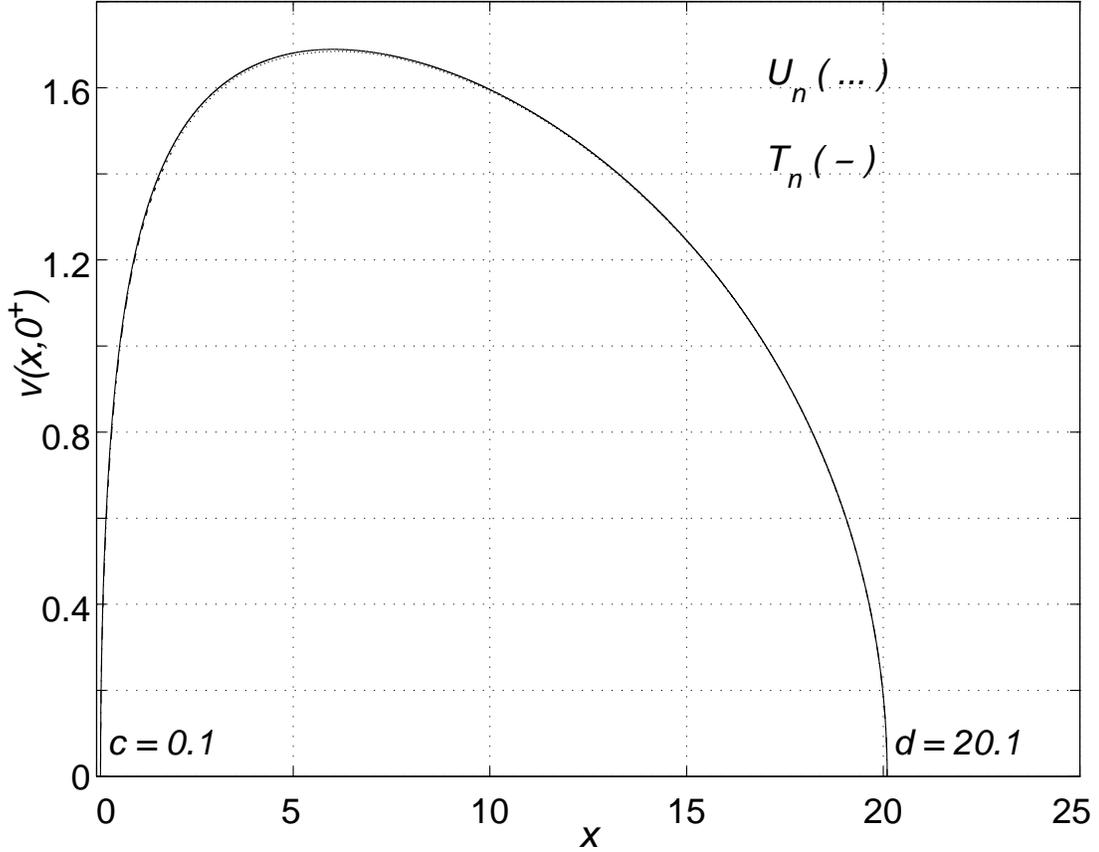}}
 \caption{Displacement profiles for a mode I crack in an infinite
 strip obtained by means of $U_n$ and $T_n$ representations
 ($N+1=15$).  Here $c = 0.1$, $d = 20.1$, $2a = 20$, and $(c+d)/(d-c)
 = 1.01$.  The crack is tilted to the left because of the ``edge effect''.}
 \label{fig2}
\end{figure}

The case $c>0$ represents an internal crack, which is the case of
interest in this work.  Based on the dominant behavior of the singular
kernels of the integral equation (\ref{intKaya2}), the solution takes
the form
\[
D(s) = R(s)\sqrt{1-s^2}\ .
\]
Here the representation function $R(s)$ is approximated in terms of
Tchebyshev polynomials of 1st and 2nd kinds, {\it i.e.}
\[
R(s) = \sum_{n=0}^Na_nU_n(s)\quad \mbox{and}
\quad R(s) = \sum_{n=0}^Nb_nT_n(s)\ .
\]
The unknown coefficients $a_n$ and $b_n$ are determined by selecting
an appropriate set of collocation points
\[
r_j = \cos\left(\frac{(2n-1)\pi}{2(N+1)}\right) \ , \ j =
1,2,\cdots,N+1; \ \ \mbox{for} \ U_n \ \mbox{representation}\ .
\]
\[
r_j = \cos\left(\frac{n\pi}{N+2}\right) \ , \ j =
1,2,\cdots,N+1; \ \ \mbox{for} \ T_n \ \mbox{representation}\ .
\]

\ 

Once the solution is obtained, the SIFs can be calculated from\footnote{Kaya and
Erdogan~\cite{KayaErdo87} do not consider the factor $\sqrt{\pi}$ in
the definition of SIFs, equations (\ref{SIFdispL}) and
(\ref{SIFdispR}).  Note that this does not affect the normalized SIFs
({\it e.g.} see \textsf{Table~\ref{tab2}}).}
\begin{eqnarray}
K_{I}(c)
 & = & \lim_{x\rightarrow c^-}\sqrt{2\pi(c-x)}\sigma_{yy}(x,0) \ , \ \ ( x
 < c ) \nonumber \\
 & = & \left(\frac{2\mu}{1+\kappa}\right)\lim_{x\rightarrow
 c^+}\frac{D(x)}{\sqrt{2\pi(x-c)}} \ , \ \ ( x > c ) \nonumber \\
 & = & \left(\frac{2\mu}{1+\kappa}\right)\sqrt{\frac{d-c}{2\pi}}R(-1)
\label{SIFdispL}
\end{eqnarray}
and
\begin{eqnarray}
K_{I}(d)
 & = & \lim_{x\rightarrow d^+}\sqrt{2\pi(x-d)}\sigma_{yy}(x,0) \ , \ \ ( x
 > d ) \nonumber \\
 & = & \left(\frac{2\mu}{1+\kappa}\right)\lim_{x\rightarrow
 d^-}\frac{D(x)}{\sqrt{2\pi(d-x)}} \ , \ \ ( x < d ) \nonumber \\
 & = & \left(\frac{2\mu}{1+\kappa}\right)\sqrt{\frac{d-c}{2\pi}}R(+1)
\label{SIFdispR}
\end{eqnarray}
which are obtained from equation (\ref{intKaya1}) by observing that
its left-hand-side gives the stress component $\sigma_{yy}(x,0)$
outside the crack interval $(c, d)$.

\ 

\begin{table}[htbp]
\centering
 \caption{Normalized stress intensity factors (SIFs) for an internal
 crack in a half-plane.  $N+1$ terms are used in approximating the
 primary variable.}
\begin{tabular}{||c|c||c|c||c|c||c|c||}      \hline
 &  & \multicolumn{2}{c||}{$U_n$ Representation} &
\multicolumn{2}{c||}{$T_n$ Representation} &
\multicolumn{2}{c||}{Kaya and Erdogan~\cite{KayaErdo87}} \\ \cline{3-8}
$\frac{d+c}{d-c}$ & $N+1$ & $\frac{K_I(c)}{p_0\sqrt{\pi(d-c)/2}}$ &
$\frac{K_I(d)}{p_0\sqrt{\pi(d-c)/2}}$ &
$\frac{K_I(c)}{p_0\sqrt{\pi(d-c)/2}}$ &
$\frac{K_I(d)}{p_0\sqrt{\pi(d-c)/2}}$ &
$\frac{K_I(c)}{p_0\sqrt{\pi(d-c)/2}}$ &
$\frac{K_I(d)}{p_0\sqrt{\pi(d-c)/2}}$ \\ \hline
1.01 & 15 & 3.6437 & 1.3292 & 3.8037 & 1.3313 & 3.6387 & 1.3298 \\ \hline
1.05 & 10 & 2.1541 & 1.2535 & 2.1920 & 1.2543 & 2.1547 & 1.2536 \\ \hline
1.1  & 10 & 1.7583 & 1.2108 & 1.7655 & 1.2111 & 1.7587 & 1.2108 \\ \hline
1.2  & 6  & 1.4637 & 1.1625 & 1.4728 & 1.1632 & 1.4637 & 1.1626 \\ \hline
1.3  & 6  & 1.3316 & 1.1331 & 1.3346 & 1.1335 & 1.3316 & 1.1331 \\ \hline
1.4  & 6  & 1.2544 & 1.1123 & 1.2556 & 1.1125 & 1.2544 & 1.1123 \\ \hline
1.5  & 4  & 1.2036 & 1.0966 & 1.2066 & 1.0969 & 1.2035 & 1.0967 \\ \hline
2.0  & 4  & 1.0913 & 1.0539 & 1.0916 & 1.0540 & 1.0913 & 1.0539 \\ \hline
3.0  & 4  & 1.0345 & 1.0246 & 1.0346 & 1.0246 & 1.0345 & 1.0246 \\ \hline
4.0  & 4  & 1.0182 & 1.0141 & 1.0182 & 1.0141 & 1.0182 & 1.0141 \\ \hline
5.0  & 4  & 1.0112 & 1.0092 & 1.0112 & 1.0092 & 1.0112 & 1.0092 \\ \hline
10.0 & 4  & 1.0026 & 1.0024 & 1.0026 & 1.0024 & 1.0026 & 1.0024 \\ \hline
20.0 & 4  & 1.0006 & 1.0006 & 1.0006 & 1.0006 & 1.0006 & 1.0006 \\ \hline
\end{tabular}
 \label{tab2}
\end{table}

\begin{figure}[ht]
\centerline{\epsfxsize=14.72 cm  \epsffile{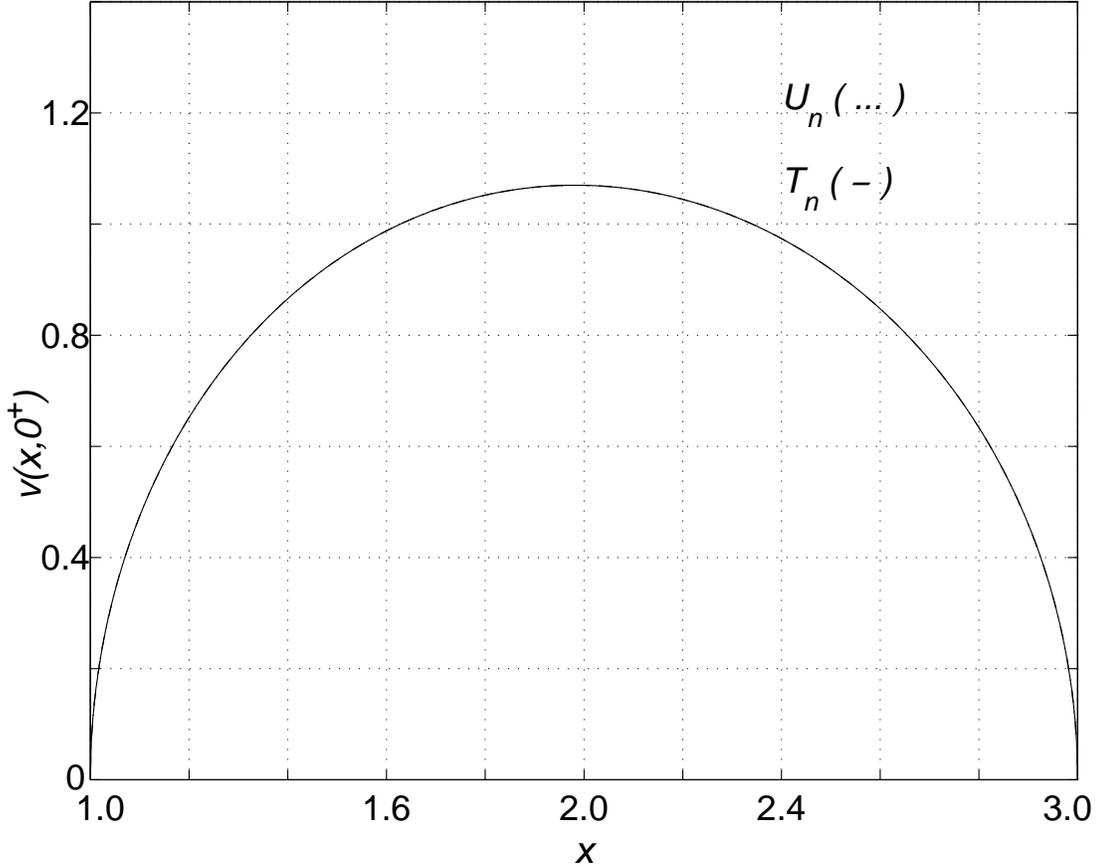}}
 \caption{Displacement profiles for a mode I crack in an infinite
 strip obtained by means of $U_n$ and $T_n$ representations
 ($N+1=8$).  Here $c = 1$, $d = 3$, $2a = 2$, and $(c+d)/(d-c)
 = 2$.}
 \label{fig3}
\end{figure}

\textsf{Table~\ref{tab2}} presents the SIFs at both tips of an
internal crack in a half-plane ($h \rightarrow \infty$) under uniform
load ($p(x) = p_0$) obtained with both $U_n$ and $T_n$
representations.  First, it is worth noting that the present SIF
results for the $U_n$ representation compare well with those reported
in Table 1 (page 114) of the paper by Kaya and
Erdogan~\cite{KayaErdo87} for the entire range of values describing
the relative position of the crack, {\it{i.e.}} $1.01 < (d+c)/(d-c) <
20$.  
Next, comparing the SIFs obtained with the
$U_n$ and $T_n$ representations in \textsf{Table~\ref{tab2}}, we note
that the results compare quite well, except when $(d+c)/(d-c) \approx
1.0$, and the discrepancy is bigger at the left-hand-side (LHS) than at the
right-hand-side (RHS) crack tip.  This occurs because of the ``edge
effect''~\cite{PaulSaifMukh93}.  If
42 terms ({\it{i.e.}} $N+1=42$) and $T_n$ representation are considered
for the case $(d+c)/(d-c) = 1.01$, then the normalized SIFs at the LHS
and RHS crack tips are 3.6437 and 1.3302, respectively.  Thus,
when there is an ``edge effect'', the results are sensitive to the
discretization adopted.  Moreover, for the same number of collocation
points, the level of accuracy attained with the $U_n$ representation
is slightly different from that with the $T_n$ representation.

\ 

\textsf{Figures~\ref{fig2}} and \textsf{~\ref{fig3}} compare the
crack profiles for $U_n$ and $T_n$ representations.  One may observe
that the displacement profiles obtained from both
representations practically agree within plotting accuracy, especially
in \textsf{Figure~\ref{fig3}}.  Note that the displacement profile in
\textsf{Figure~\ref{fig2}} is tilted to the left because of the ``edge
effect''.  Such effect is negligible in \textsf{Figure~\ref{fig3}}.

\subsection{\normalsize{Mode III Crack Problem in
Nonhomogeneous Materials~\cite{Erdogan85}}}

\begin{figure}[ht]
\centerline{\epsfxsize=14.72 cm  \epsffile{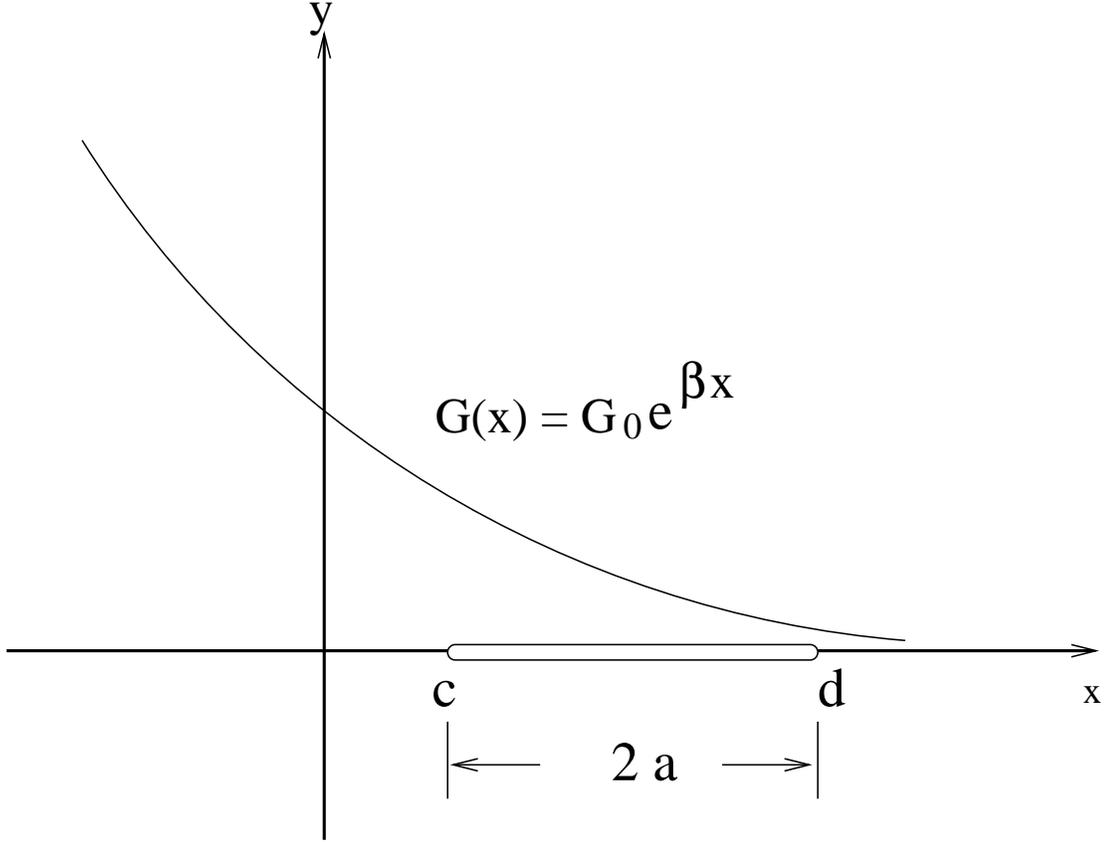}}
 \caption{The half plane of the antiplane shear problem for
 nonhomogeneous material with shear modulus $G(x)=G_0\myE^{\beta x}$.}
 \label{fig4}
\end{figure}

Consider the antiplane shear problem for the nonhomogeneous material
shown in \textsf{Figure~\ref{fig4}}
with shear modulus variation given by
\begin{equation}
G(x) = G_0\myE^{\beta x} \ ,
\end{equation}
where $G_0$ and $\beta$ are material constants.
Erdogan~\cite{Erdogan85} has studied this  problem in order to
investigate the singular nature of the crack-tip stress field
in bonded nonhomogeneous materials under antiplane shear loading.
He uses a slope formulation\footnote{The governing
integral equations are described by relationships (20), (21), and (22)
of Reference~\cite{Erdogan85}.}, while here we use a
displacement formulation.
To understand what can be gained through the displacement
formulation, we first state the governing {\sf PDE} and the boundary
conditions for the crack problem:
\begin{equation}
\begin{array}{ll}
\nabla^{2}w(x,y) + \beta\frac{\partial w(x,y)}{\partial x} = 0\ ,
& \quad -\infty<x<\infty\ , \quad y\geq 0\ , \\
w(x,0) = 0\ , & \quad x\notin [c, d]\ , \\
\sigma_{yz}(x, 0^+) = p(x) \ , & \quad x\in (c, d)\ ,
\end{array}
\end{equation}
where $p(x)$ is the traction function along the crack surfaces $(c,
d)$; and because of symmetry, only the upper half plane $y > 0$ is considered.
By the Fourier transform we write
$w(x,y)$ as
\begin{equation}
w(x,y) =
\frac{1}{\sqrt{2\pi}}\int_{-\infty}^{\infty}
\left[A(\xi)\myE^{\lambda(\xi)y}\right]\myE^{ix\xi}d\xi\ ,
\label{functionA}
\end{equation}
where $A(\xi)$ is to be determined by the boundary conditions, and
\begin{equation}
[\lambda(\xi)]^2 = \xi^2 + i\beta\xi\ .
\end{equation}
Because of the far field boundary condition, $\lim_{y
\rightarrow\infty}w(x,y) = 0$, $\lambda(\xi)$ is found to have a
non-negative real part which can be expressed as
\begin{equation}
\lambda(\xi) = \frac{-1}{\sqrt{2}}\sqrt{\sqrt{\xi^4+\beta\xi^2}+\xi^2}
-
\frac{i}{\sqrt{2}}{\mathsf S}(\beta\xi)\sqrt{\sqrt{\xi^4+\beta\xi^2}-\xi^2}\ ,
\end{equation}
where the sign function ${\mathsf S}(\cdot)$ is defined as
\begin{equation}
{\mathsf S}(\eta) = \left\{
   \begin{array}{rl}
 1\ ,  &\quad \eta > 0  \\
 0\ ,  &\quad \eta = 0  \\
-1\ ,  &\quad \eta < 0\ .
   \end{array}
                              \right.
\label{signOfBeta}
\end{equation}
By applying the inverse Fourier transform to equation
(\ref{functionA}), one finds 
\begin{equation}
A(\xi) =
\frac{1}{\sqrt{2\pi}}\int_{c}^{d}w(t,0)\myE^{it\xi}dt\ ,
\end{equation}
which leads to the following integral equation\footnote{Note that
equations (\ref{tradeOff1}) and (\ref{tradeOff2})
correspond to equations (20) and (21) in reference~\cite{Erdogan85},
respectively.  However, the present notation is different
from~\cite{Erdogan85}, in which the dummy variable used for the
Fourier transform is $\alpha$, $(a, b)$ stands for the crack surfaces, and $m(\alpha)$
corresponds to our $\lambda(\xi)$.}:
\begin{equation}
\sigma_{yz}(x, 0) = p(x) =
\frac{G(x)}{2\pi}\int_{c}^{d}k(x,t)w(t,0)dt\ ,
\label{tradeOff1}
\end{equation}
where
\begin{equation}
k(x,t) = \lim_{y\rightarrow0^+}\int_{-\infty}^{\infty}
\left[\lambda(\xi)\myE^{\lambda(\xi)y}\right]\myE^{i(t-x)\xi}d\xi\ .
\label{tradeOff2}
\end{equation}

\ 

The trade-off between a displacement-based and a slope-based
formulation can be seen here if one recalls the issue
regarding decomposition of the $K(\xi)$ described in
(\ref{decompose}).
In summary,
\begin{equation}
K(\xi) =
\lim_{y\rightarrow0^+}\frac{\lambda(\xi)}{\xi}\myE^{\lambda(\xi)y}
= \frac{\lambda(\xi)}{\xi}\ , \quad \mbox{if slope formulation is used,}
\label{tradeOff3}
\end{equation}
and
\begin{equation}
K(\xi) = \lim_{y\rightarrow0^+}\lambda(\xi)\myE^{\lambda(\xi)y} =
\lambda(\xi)\ , \quad \mbox{if displacement formulation is adopted.}
\label{tradeOff4}
\end{equation}
The decomposition of $K(\xi)$ in (\ref{tradeOff3})
is difficult because of the term $\xi$ in the denominator.  On the
other hand, the decomposition of $K(\xi)$ in (\ref{tradeOff4}) can be
achieved through a simple asymptotic analysis:
\begin{equation}
\mbox{\textsf{R}eal part of }\lambda(\xi)
= \frac{-1}{\sqrt{2}}\sqrt{\sqrt{\xi^4+\beta\xi^2}+\xi^2}
\stackrel{|\xi|\rightarrow\infty}{\ \ \longrightarrow\ \ }-|\xi|\ ,
\end{equation}

\begin{equation}
i \times
\mbox{\textsf{I}maginary part of }\lambda(\xi) = 
\frac{-i}{\sqrt{2}}{\mathsf S}(\beta\xi)\sqrt{\sqrt{\xi^4+\beta\xi^2}-\xi^2}
\stackrel{|\xi|\rightarrow\infty}{\ \ \longrightarrow\ \ }-\frac{i\beta}{2}\frac{|\xi|}{\xi}\ .
\end{equation}
Thus the governing hypersingular integral equation is found to be
\begin{equation}
\frac{G(x)}{2\pi}\int_{c}^{d}\left[\frac{2}{(t-x)^2}+
\frac{\beta}{t-x}+N(x,t)\right]D(t)dt\ =\ p(x)\ ,\quad c<x<d\ ,
\label{IntEqnExp2}
\end{equation}
where we have let
\begin{equation}
D(t)\ =\ w(t,0)\ ,
\end{equation}
and the nonsingular kernel is
\begin{eqnarray}
N(x,t) & = &
\int_{0}^{\infty}\left\{
\frac{-\beta^2\sqrt{\xi}\cos[(t-x)\xi]}{\left(\sqrt{\xi}+\frac{1}{\sqrt{2}}\sqrt{\sqrt{\xi^2+\beta^2}+\xi}\right)
\left(\xi+\sqrt{\xi^2+\beta^2}
\right)} \ + \right.\nonumber \\
& \ & \left.\frac{-\beta^4\sin[(t-x)\xi]}{\left(\beta+\sqrt{2}{\mathsf S}(\beta)\sqrt{\sqrt{\xi^4+\beta^2\xi^2}-\xi^2}\right)
\left(2\xi^2+\beta^2+2\sqrt{\xi^4+\beta^2\xi^2}
\right)}\right\}d\xi \ .
\label{RegKrlExp2}
\end{eqnarray}
Recall that the function ${\mathsf S}(\cdot)$ is defined by
equation (\ref{signOfBeta}).  As a consistency check, note that if
$\beta = 0$, then both the Cauchy
singular kernel $\beta / (x-t)$ and the nonsingular kernel $N(x,t)$
will be dropped from equation (\ref{IntEqnExp2}) so that equation
(\ref{QuadraForm}) is recovered.

\begin{figure}[ht]
\centerline{\epsfxsize= 12.72 cm \epsffile{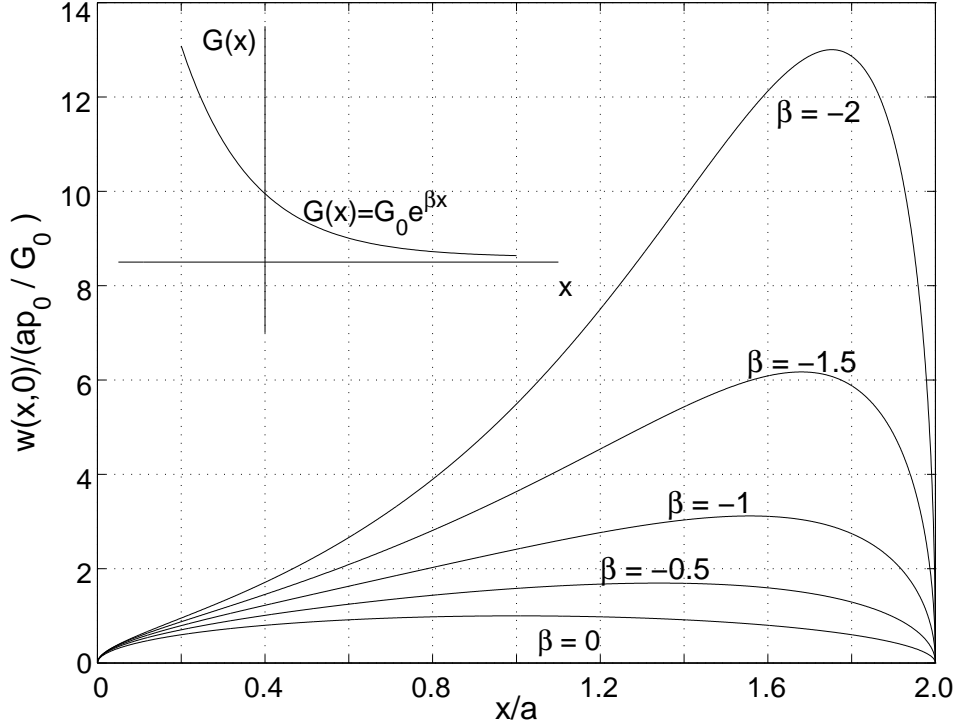}}
 \caption{The half plane of the antiplane shear problem for
 nonhomogeneous material with shear modulus $G(x)=G_0\myE^{\beta x}$.
 The cracks are tilted to the right because of material nonhomogeneity.}
 \label{fig5}
\end{figure}

{\sf Figure~\ref{fig5}} shows numerical results for displacement
profiles considering a crack with uniformly applied shear tractions
$\sigma_{yz}(x,0) = -p_0 \ (|x|<a)$, and various values of the
material parameter $\beta$.
Note that the cracks are tilted to the right because of material
nonhomogeneity.  Further numerical results, including SIFs at both
tips of the crack and corresponding displacement profiles, are given
by Chan {\it et. al.}~\cite{ChanPaulFann99a}.  From a numerical point
of view, they have shown that essentially the same results are
obtained either by $U_n$ or $T_n$ representations
for this specific problem~\cite{ChanPaulFann99a}.

\subsection{\normalsize{Gradient Elasticity Applied to
Mode III Cracks~\cite{PaulChanFann99a}}}
One of the most relevant aspects of the formulas derived in
{\textsf{Sections 4}} and {\textsf 5} is the evaluation of hypersingular
integrals such as $I_{\alpha}(T_n,m,r)$ for $m \geq 2$ in the weight
function $W(s)$ given by equations (\ref{repDnsty}) and
(\ref{repWt}).  This example illustrates this point for the case
$m=2$.

\

Paulino {\it et. al.}~\cite{PaulChanFann99a} have presented a
hypersingular integral equation formulation for a mode III crack in a
material described by constitutive equations of gradient elasticity
with both volumetric and surface energy gradient dependent terms.  A
similar study, using a different approach, has been conducted by
Vardoulakis {\it et. al.}~\cite{VardExadAifa96} and Exadaktylos {\it
et. al.}~\cite{ExadVardAifa96}.  For this problem, the governing {\sf
PDE} is
\begin{equation}
-\ell^2 \nabla^{4} w + \nabla^{2} w = 0 \ ,
\label{pdeIII}
\end{equation}
where $\ell$ is the characteristic length of the
material associated to volumetric strain-gradient terms, and $w$ is
the antiplane shear displacement.  The boundary conditions are
\begin{equation}
 \begin{array}{ll}
 \sigma_{yz}(x,0) = p(x)  \  , &  |x| < a \\
 \mu_{yyz}(x,0) =  0  \ , & -\infty < x < \infty  \\
 w(x,0) =  0  \ , & |x| > a  \ , \\
 \end{array}
\label{bdsIII}
\end{equation}
where the notation of references~\cite{ExadVardAifa96,PaulChanFann99a,
VardExadAifa96} is adopted here.
Enforcing the governing equation (\ref{pdeIII}), imposing the boundary
conditions (\ref{bdsIII}), taking account of symmetry along the
$x$-axis, and using Fourier transform method, Paulino {\it
et. al.}~\cite{PaulChanFann99a} have obtained the following governing Fredholm
hypersingular integral equation
\begin{equation}
\frac{G}{\pi}\int_{-a}^{a}\left\{ \frac{-2\ell^2}{(t-x)^3} +
\frac{1-(\frac{\myLp}{2\ell})^2}{t-x} + k(x,t)\right\} \phi (t) \ dt +
\frac{G\ell '}{2} \phi '(x)
=  p(x), \ |x|<a 
\label{intEqnIII}
\end{equation}
with singularity $\alpha = 3$,
where the slope function $\phi(x)$ is defined to be
\begin{equation}
\phi(x) = \frac{\partial w(x,0)}{\partial x}
\label{exp3dnsty}
\end{equation}
which satisfies the single-valuedness condition
\begin{equation}
\int_{-a}^{a} \phi(t)dt = 0 \ ,
\end{equation}
for the solution of the fracture mechanics problem.  Here, $k(x,t)$ is
the nonsingular kernel given by
\begin{equation}
k(x,t)=\int_{0}^{\infty}\frac{\frac{\myLp}{2}\xi \left(\sqrt{\frac{\ell^{2}\xi^2
+1}{\ell^2}}-\xi
\right)-\frac{1}{4}(\frac{\myLp}{\ell})^2\left(\sqrt{\frac{\ell^{2}\xi^2
+1}{\ell^2}}-\xi
\right)+\frac{1}{4}\frac{(\myLp)^3}{\ell^4}}{\frac{\myLp}{\ell^2}-\left(\sqrt{\frac{\ell^{2}\xi^2
+1}{\ell^2}}+\xi \right)}\sin[\xi (t-x)]d\xi \ ;
\label{nSnglr3}
\end{equation}
where $(-a,\  a)$ stands for the crack surfaces; $\myLp$ is
the characteristic length responsible for surface strain-gradient
terms; $G$ is
the shear modulus of the material; $p(x)$ is
the known loading function; $t$ is the integration variable, and $x$
is the collocation variable.

\ 

The behavior of the solution in terms of the density function $\phi(t)$,
can be expressed as
\begin{equation}
\phi(t)\ = \ R(t)(a^2-t^2)^{\frac{3}{2}} \ ,
\label{exp3Dnsty}
\end{equation}
where $R(t)$ is taken to be an expansion of Tchebyshev
polynomials of first kind ($T_n$).  This example motivates the
whole work of this paper because the analytical evaluation of the
hypersingular integral $I_{3}(T_n,2,r)$ is needed for successfully
solving the governing integral equation (\ref{intEqnIII}).  The
expression for $I_{3}(T_n,2,r)$ is given by (\ref{hyperCube}).

\

Note that the unknown density function is taken to be the first
derivative of the displacement function, as described by equation
(\ref{exp3dnsty}).  For this particular example, the decomposition of
the original kernel into singular and nonsingular parts, stated by
equation (\ref{decompose}), can be accomplished by means of partial
fractions~\cite{PaulChanFann99a}.
In general, this step of
decomposition is not an easy task, as discussed in {\textsf{Example 2}}.
An alternative is to consider the
displacement function $w(x,0)$ as the unknown density function.
In this case, a hypersingular integral equation
with $\alpha = 4$ is obtained.  With order of singularity $\alpha =
4$, the behavior of the density function $D(t)$ in
equation (\ref{int_xt2}) can be expressed by~\cite{PaulChanFann99a}
\begin{equation}
D(t)\ = \ R(t)(a^2-t^2)^{\frac{5}{2}} \ .
\label{exp3Dnsty2}
\end{equation}
Thus one needs to evaluate the hypersingular integral
$I_{4}(T_n,3,r)$ in order to implement the numerical approximation.
The expression for $I_{4}(T_n,3,r)$ is given by equation (\ref{myRefExp3}).

\

\begin{table}[htbp]
\centering
 \caption{Stress intensity factors
 $K_{III}(a)=3\sqrt{\pi a}(\frac{\ell}{a})^2G\sum_{n=0}^{N}a_n$,
 where $a_n$ are the coefficients of the $T_n$ expantion to $R(t)$ in
 equation (\ref{exp3Dnsty}).  Here $\myLp = 0$}
\begin{tabular}{||c|c|c|c|c|c|c|c||}      \hline
N & $\ell$=0.8 & $\ell$=0.5 & $\ell$=0.2 & $\ell$=0.1 & $\ell$=0.05
& $\ell$=0.01 & $\ell$=0.005  \\ \hline
11 & 20.3131 & 15.8292 & 7.4396 & 4.5116 & 2.6342 & 0.1282 & 0.0319  \\ \hline

21 & 11.8757 & 9.5632 & 4.4791 & 2.1538 & 0.9541 & 0.0898 & 0.1602  \\ \hline

31 & 11.6607 & 9.3937 & 4.3878 & 2.0856 & 0.9204 & 0.1649 & 0.0404  \\ \hline

41 & 11.6665 & 9.3983 & 4.3902 & 2.0878 & 0.9246 & 0.1378 & 0.0682  \\ \hline

51 & 11.6667 & 9.3984 & 4.3903 & 2.0878 & 0.9247 & 0.1399 & 0.0658  \\ \hline

61 & 11.6667 & 9.3984 & 4.3903 & 2.0878 & 0.9247 & 0.1400 & 0.0653  \\ \hline

71 & 11.6667 & 9.3984 & 4.3903 & 2.0878 & 0.9247 & 0.1399 & 0.0654  \\ \hline

81 & 11.6667 & 9.3984 & 4.3903 & 2.0878 & 0.9247 & 0.1399 & 0.0654  \\ \hline

91 & 11.6667 & 9.3984 & 4.3903 & 2.0878 & 0.9247 & 0.1399 & 0.0654  \\ \hline

101 & 11.6667 & 9.3984 & 4.3903 & 2.0878 & 0.9247 & 0.1399 & 0.0654  \\ \hline

\end{tabular}
 \label{tab4}
\end{table}

\begin{figure}[ht]
\centerline{\epsfxsize=14.72 cm \epsffile{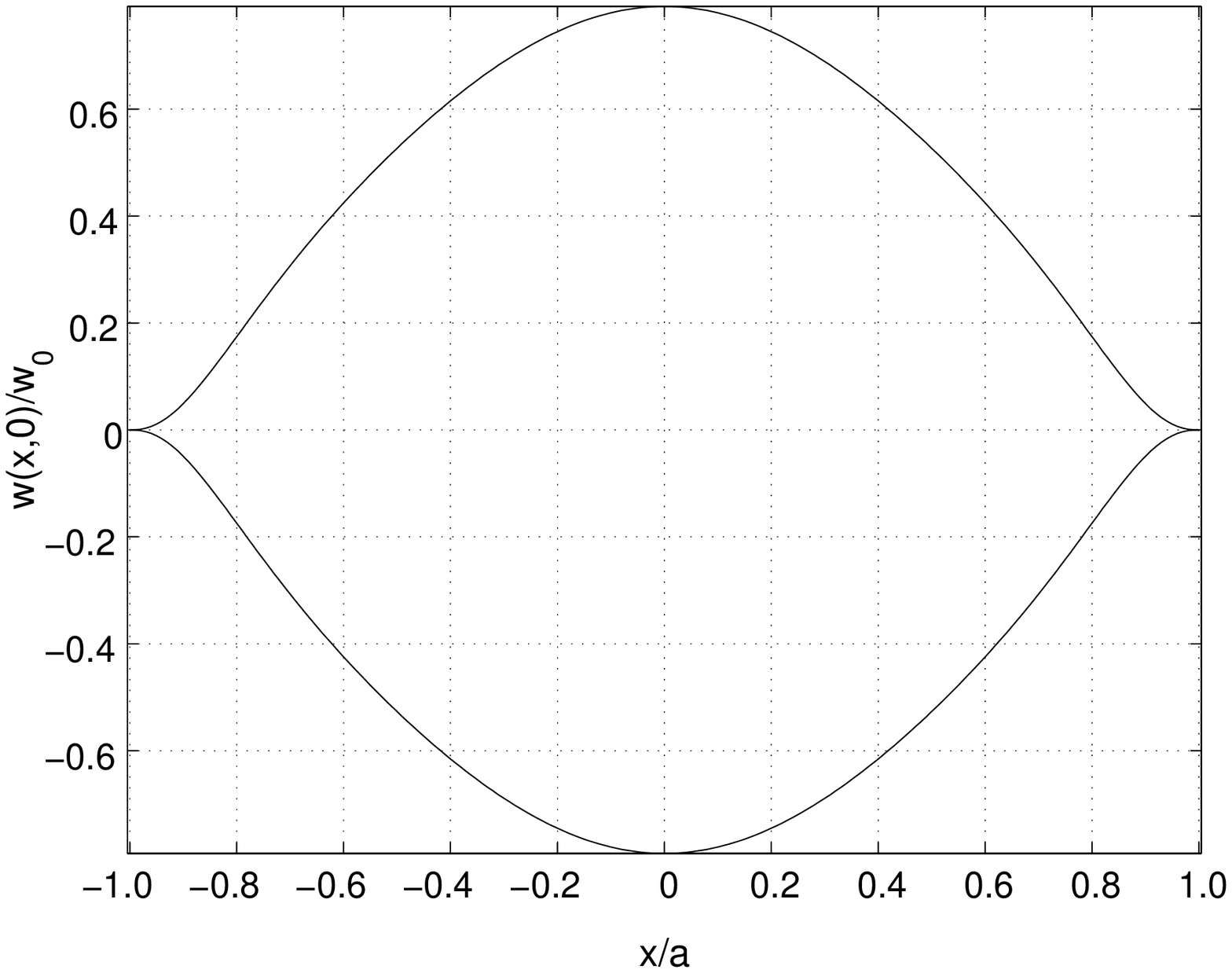}}
 \caption{Full crack displacement profile at $\ell/a = 0.2$ and $\myLp/a =
 0.1$ under uniform crack surface antiplane loading
 $\sigma_{yz}(x,0)=-\sigma_0$; $w_0 = a\sigma_0/G_0$.}
 \label{fig6}
\end{figure}

The numerical results are presented in terms of SIFs ({\sf Table~\ref{tab4}})
and displacement profiles ({\sf Figure~\ref{fig6}}) by considering the
slope function $\phi(x)$ as unknown in equation (\ref{intEqnIII})
and the $T_n$ expansion to $R(t)$
in equation (\ref{exp3Dnsty}).  {\sf Table~\ref{tab4}} shows the
convergence of the SIFs by choosing different values of volumetric
gradient dependent term, $\ell$, and letting surface energy gradient dependent
term $\myLp = 0$.  The displacement profile in {\sf Figure~\ref{fig6}}
shows the ``cusping'' phenomenon at the crack tips, which is also described by 
Barenblatt's ``cohesive zone'' theory~\cite{Barenb62}.  The interesting point is that
the cusping crack obtained here is a natural outcome of the higher order
gradient elasticity theory.
Further numerical results and discussions are provided in
reference~\cite{PaulChanFann99a}.


\section{\large{Concluding Remarks}}
Closed form analytical solutions are provided here for a broad class
of improper integrals with hypersingular kernels and
density functions approximated
by means of Tchebyshev polynomials.  Whenever possible, the symbolical and
numerical tools of the
computer algebra software MAPLE
\footnote{Maple can evaluate some CPV integrals.  However, in general,
computer algebra systems are very limited with respect to hypersingular
integrals.}~\cite{Heck96, Waterloo95a, Waterloo95b} have
been used to verify the proposed
solutions. A systematic approach for evaluating
integrals when higher order singularities is also given in the present paper.

\

The examples involve crack problems and aspects such as LEFM,
nonhomogeneous materials, and gradient elasticity theory (see {\textsf{
Table~\ref{tab1}}}).  All these problems
are solved by means of Fredholm hypersingular integral equation formulations.
When classical elasticity is used, both $T_n$ and $U_n$ representations lead
to essentially the same numerical results.  For a crack problem in
nonhomogeneous material, the difficulty that arises in splitting the
singular and nonsingular parts from the original kernels can be circumvented
by means of displacement-based, rather than slope-based,
formulation.

\ 

As a closing remark, we note that as material property variation in space
and higher order gradient continuum theories are considered, the formulation
of the crack problem and the associated kernels become quite involved.
Thus, better analytical and numerical techniques are needed to successfully
solve the governing singular integral equations.  This paper is a combination
in this sense.


\ 

\newpage

\vspace{1.0in}

\appendix

{\flushleft \LARGE{\textsf{Appendix}}}
\section{Integrals associated with lower order $n$}
The general formula given in the text, {\it e.g.} equations
(\ref{boxI}), (\ref{boxII}), (\ref{boxIII}), (\ref{boxIV}), (\ref{boxV}),
(\ref{boxVI}), (\ref{boxVII}), and (\ref{boxVIII}) are only valid
above certain values of $n$.  Thus
the goal of this {\textsf{Appendix}} is to provide the expressions for
integrals associated with lower order $n$.

\subsection{$I_{1}(T_n,1,r), \ n=0, \ 1$}
\begin{equation}
\frac{1}{\pi}\int_{-1}^{1}\cpv\frac{T_0(s)\sqrt{1-s^2}}{s-r}ds = 
-r \ , \ \ |r|<1
\end{equation}

\begin{equation}
\frac{1}{\pi}\int_{-1}^{1}\cpv\frac{T_1(s)\sqrt{1-s^2}}{s-r}ds = 
-r^2 + \frac{1}{2} \ , \ \ |r|<1
\end{equation}

\subsection{$I_{1}(T_n,2,r), \ n=0\cdots 5$}

\begin{equation}
\frac{1}{\pi}\int_{-1}^{1}\cpv\frac{T_0(s)(1-s^2)^{\frac{3}{2}}}{s-r}ds = r^3
-\frac{3}{2}r \ , \ \ |r|<1
\label{cauZero}
\end{equation}

\begin{equation}
\frac{1}{\pi}\int_{-1}^{1}\cpv\frac{T_1(s)(1-s^2)^{\frac{3}{2}}}{s-r}ds =
r^4 - \frac{3}{2}r^2 +\frac{3}{8} \ , \ \ |r|<1
\label{cauOne}
\end{equation}

\begin{equation}
\frac{1}{\pi}\int_{-1}^{1}\cpv\frac{T_2(s)(1-s^2)^{\frac{3}{2}}}{s-r}ds =
2r^5 - 4r^3 +\frac{9}{4}r \ , \ \ |r|<1
\label{cauTwo}
\end{equation}

\begin{equation}
\frac{1}{\pi}\int_{-1}^{1}\cpv\frac{T_3(s)(1-s^2)^{\frac{3}{2}}}{s-r}ds =
4r^6 - 9r^4 +6r^2 -\frac{7}{8} \ , \ \ |r|<1
\label{cauThree}
\end{equation}

\begin{equation}
\frac{1}{\pi}\int_{-1}^{1}\cpv\frac{T_4(s)(1-s^2)^{\frac{3}{2}}}{s-r}ds =
8r^7 - 20r^5 + 16r^3 - 4r \ , \ \ |r|<1
\label{cauFour}
\end{equation}

\begin{equation}
\frac{1}{\pi}\int_{-1}^{1}\cpv\frac{T_5(s)(1-s^2)^{\frac{3}{2}}}{s-r}ds =
16r^8 - 44r^6 + 41r^4 -14r^2 + 1 \ , \ \ |r|<1
\label{cauFive}
\end{equation}

\subsection{$I_{1}(T_n,3,r), \ n=0\cdots 5$}
\begin{equation}
\frac{1}{\pi}\int_{-1}^{1}\cpv\frac{T_0(s)(1-s^2)^{\frac{5}{2}}}{s-r}ds =
- \frac{15}{8}r + \frac{5}{2}r^3  - r^5
\end{equation}

\begin{equation}
\frac{1}{\pi}\int_{-1}^{1}\cpv\frac{T_1(s)(1-s^2)^{\frac{5}{2}}}{s-r}ds =
\frac{5}{16} - \frac{15}{8}r^2 + \frac{5}{2}r^4  - r^6
\end{equation}

\begin{equation}
\frac{1}{\pi}\int_{-1}^{1}\cpv\frac{T_2(s)(1-s^2)^{\frac{5}{2}}}{s-r}ds =
\frac{5}{12}r - \frac{25}{4}r^3 + 6r^5 - 2r^7
\end{equation}

\begin{equation}
\frac{1}{\pi}\int_{-1}^{1}\cpv\frac{T_3(s)(1-s^2)^{\frac{5}{2}}}{s-r}ds =
-\frac{25}{32} + \frac{55}{8}r^2 - 15r^4  + 13r^6 - 4r^8
\end{equation}

\begin{equation}
\frac{1}{\pi}\int_{-1}^{1}\cpv\frac{T_4(s)(1-s^2)^{\frac{5}{2}}}{s-r}ds =
-\frac{65}{16}r +20r^3 - 36r^5 + 28r^7 -8r^9
\end{equation}

\begin{equation}
\frac{1}{\pi}\int_{-1}^{1}\cpv\frac{T_5(s)(1-s^2)^{\frac{5}{2}}}{s-r}ds =
\frac{31}{32} -15r^2 + 55r^4  - 85r^6 + 60r^8 -16r^{10}
\end{equation}

\subsection{$I_{1}(U_n,3,r), \ n=0\cdots 5$}
\begin{equation}
\frac{1}{\pi}\int_{-1}^{1}\cpv\frac{U_0(s)(1-s^2)^{\frac{5}{2}}}{s-r}ds =
- \frac{15}{8}r + \frac{5}{2}r^3  - r^5
\end{equation}

\begin{equation}
\frac{1}{\pi}\int_{-1}^{1}\cpv\frac{U_1(s)(1-s^2)^{\frac{5}{2}}}{s-r}ds =
\frac{5}{8} - \frac{15}{4}r^2 + 5r^4  - 2r^6
\end{equation}

\begin{equation}
\frac{1}{\pi}\int_{-1}^{1}\cpv\frac{U_2(s)(1-s^2)^{\frac{5}{2}}}{s-r}ds =
\frac{25}{8}r - 10r^3 + 11r^5 - 4r^7
\end{equation}

\begin{equation}
\frac{1}{\pi}\int_{-1}^{1}\cpv\frac{U_3(s)(1-s^2)^{\frac{5}{2}}}{s-r}ds =
-\frac{15}{16} + 10r^2 - 25r^4  + 24r^6 - 8r^8
\end{equation}

\begin{equation}
\frac{1}{\pi}\int_{-1}^{1}\cpv\frac{U_4(s)(1-s^2)^{\frac{5}{2}}}{s-r}ds =
-5r +30r^3 - 61r^5 + 52r^7 -16r^9
\end{equation}

\begin{equation}
\frac{1}{\pi}\int_{-1}^{1}\cpv\frac{U_5(s)(1-s^2)^{\frac{5}{2}}}{s-r}ds =
1 -20r^2 + 85r^4  - 146r^6 + 112r^8 -32r^{10}
\end{equation}

\subsection{$I_{2}(T_n,3,r), \ n=0\cdots 6$}
\begin{equation}
\frac{1}{\pi}\int_{-1}^{1}\hypsng\frac{T_0(s)(1-s^2)^{\frac{5}{2}}}{(s-r)^2}ds
= -5r^4 + \frac{15}{2}r^2 -\frac{15}{8} \ , \ \ |r|<1
\end{equation}

\begin{equation}
\frac{1}{\pi}\int_{-1}^{1}\hypsng\frac{T_1(s)(1-s^2)^{\frac{5}{2}}}{(s-r)^2}ds
= -6r^5 + 10r^3 - \frac{15}{4}r \ , \ \ |r|<1
\end{equation}

\begin{equation}
\frac{1}{\pi}\int_{-1}^{1}\hypsng\frac{T_2(s)(1-s^2)^{\frac{5}{2}}}{(s-r)^2}ds
= -14r^6 + 30r^4 - \frac{75}{4}r^2 + \frac{5}{2} \ , \ \ |r|<1
\end{equation}

\begin{equation}
\frac{1}{\pi}\int_{-1}^{1}\hypsng\frac{T_3(s)(1-s^2)^{\frac{5}{2}}}{(s-r)^2}ds
= -32r^7 + 78r^5 - 60r^3 + \frac{55}{4}r \ , \ \ |r|<1
\end{equation}

\begin{equation}
\frac{1}{\pi}\int_{-1}^{1}\hypsng\frac{T_4(s)(1-s^2)^{\frac{5}{2}}}{(s-r)^2}ds
= -72r^8 + 196r^6 - 180r^4 + 60r^2 - \frac{65}{16}\ , \ \ |r|<1
\end{equation}

\begin{equation}
\frac{1}{\pi}\int_{-1}^{1}\hypsng\frac{T_5(s)(1-s^2)^{\frac{5}{2}}}{(s-r)^2}ds
= -160r^9 +480r^7 -510r^5 +220r^3 -30r \ , \ \ |r|<1
\end{equation}

\begin{equation}
\frac{1}{\pi}\int_{-1}^{1}\hypsng\frac{T_6(s)(1-s^2)^{\frac{5}{2}}}{(s-r)^2}ds
= -352r^{10} +1152r^8 -1386r^6 +730r^4 -150r^2 +6\ , \ \ |r|<1
\end{equation}

\subsection{$I_{2}(U_n,3,r), \ n=0\cdots 6$}
\begin{equation}
\frac{1}{\pi}\int_{-1}^{1}\hypsng\frac{U_0(s)(1-s^2)^{\frac{5}{2}}}{(s-r)^2}ds
= -5r^4 + \frac{15}{2}r^2 -\frac{15}{8} \ , \ \ |r|<1
\end{equation}

\begin{equation}
\frac{1}{\pi}\int_{-1}^{1}\hypsng\frac{U_1(s)(1-s^2)^{\frac{5}{2}}}{(s-r)^2}ds
= -12r^5 + 20r^3 - \frac{15}{2}r \ , \ \ |r|<1
\end{equation}

\begin{equation}
\frac{1}{\pi}\int_{-1}^{1}\hypsng\frac{U_2(s)(1-s^2)^{\frac{5}{2}}}{(s-r)^2}ds
= -28r^6 + 55r^4 - 30r^2 + \frac{25}{8} \ , \ \ |r|<1
\end{equation}

\begin{equation}
\frac{1}{\pi}\int_{-1}^{1}\hypsng\frac{U_3(s)(1-s^2)^{\frac{5}{2}}}{(s-r)^2}ds
= -64r^7 + 144r^5 - 100r^3 + 20r \ , \ \ |r|<1
\end{equation}

\begin{equation}
\frac{1}{\pi}\int_{-1}^{1}\hypsng\frac{U_4(s)(1-s^2)^{\frac{5}{2}}}{(s-r)^2}ds
= -144r^8 + 364r^6 - 305r^4 + 90r^2 - 5\ , \ \ |r|<1
\end{equation}

\begin{equation}
\frac{1}{\pi}\int_{-1}^{1}\hypsng\frac{U_5(s)(1-s^2)^{\frac{5}{2}}}{(s-r)^2}ds
= -320r^9 +896r^7 -876r^5 +340r^3 -40r \ , \ \ |r|<1
\end{equation}

\begin{equation}
\frac{1}{\pi}\int_{-1}^{1}\hypsng\frac{U_6(s)(1-s^2)^{\frac{5}{2}}}{(s-r)^2}ds
= -704r^{10} +2160r^8 -2408r^6 +1155r^4 -210r^2 +7\ , \ \ |r|<1
\end{equation}

\subsection{$I_{3}(T_n,2,r), \ n=0\cdots 4$}
\begin{equation}
\frac{1}{\pi}\int_{-1}^{1}\hypsng\frac{T_0(s)(1-s^2)^{\frac{3}{2}}}{(s-r)^3}ds
= 3r \ , \ \ |r|<1
\end{equation}

\begin{equation}
\frac{1}{\pi}\int_{-1}^{1}\hypsng\frac{T_1(s)(1-s^2)^{\frac{3}{2}}}{(s-r)^3}ds
= 6r^2 - \frac{3}{2} \ , \ \ |r|<1
\end{equation}

\begin{equation}
\frac{1}{\pi}\int_{-1}^{1}\hypsng\frac{T_2(s)(1-s^2)^{\frac{3}{2}}}{(s-r)^3}ds
= 20r^3 - 12r \ , \ \ |r|<1
\end{equation}

\begin{equation}
\frac{1}{\pi}\int_{-1}^{1}\hypsng\frac{T_3(s)(1-s^2)^{\frac{3}{2}}}{(s-r)^3}ds
= 60r^4 - 54r^2 +6 \ , \ \ |r|<1
\end{equation}

\begin{equation}
\frac{1}{\pi}\int_{-1}^{1}\hypsng\frac{T_4(s)(1-s^2)^{\frac{3}{2}}}{(s-r)^3}ds
= 168r^5 - 200r^3 +48r \ , \ \ |r|<1
\end{equation}

\subsection{$I_{3}(T_n,3,r), \ n=0\cdots 7$}
\begin{equation}
\frac{1}{\pi}\int_{-1}^{1}\hypsng\frac{T_0(s)(1-s^2)^{\frac{5}{2}}}{(s-r)^3}ds
= -10r^3 + \frac{15}{2}r \ , \ \ |r|<1
\end{equation}

\begin{equation}
\frac{1}{\pi}\int_{-1}^{1}\hypsng\frac{T_1(s)(1-s^2)^{\frac{5}{2}}}{(s-r)^3}ds
= -15r^4 + 15r^2 - \frac{15}{8} \ , \ \ |r|<1
\end{equation}

\begin{equation}
\frac{1}{\pi}\int_{-1}^{1}\hypsng\frac{T_2(s)(1-s^2)^{\frac{5}{2}}}{(s-r)^3}ds
= -42r^5 + 60r^3 - \frac{75}{4}r \ , \ \ |r|<1
\end{equation}

\begin{equation}
\frac{1}{\pi}\int_{-1}^{1}\hypsng\frac{T_3(s)(1-s^2)^{\frac{5}{2}}}{(s-r)^3}ds
= -112r^6 + 195r^4 - 90r^2 + \frac{55}{8} \ , \ \ |r|<1
\end{equation}

\begin{equation}
\frac{1}{\pi}\int_{-1}^{1}\hypsng\frac{T_4(s)(1-s^2)^{\frac{5}{2}}}{(s-r)^3}ds
= -288r^7 + 588r^5 - 360r^3 + 60r \ , \ \ |r|<1
\end{equation}

\begin{equation}
\frac{1}{\pi}\int_{-1}^{1}\hypsng\frac{T_5(s)(1-s^2)^{\frac{5}{2}}}{(s-r)^3}ds
= -720r^8 +1680r^6 -1275r^4 +330r^2 -15 \ , \ \ |r|<1
\end{equation}

\begin{equation}
\frac{1}{\pi}\int_{-1}^{1}\hypsng\frac{T_6(s)(1-s^2)^{\frac{5}{2}}}{(s-r)^3}ds
= -1760r^9 +4608r^7 -4158r^5 +1460r^3 -150r \ , \ \ |r|<1
\end{equation}

\begin{eqnarray}
\lefteqn{\frac{1}{\pi}\int_{-1}^{1}\hypsng
\frac{T_7(s)(1-s^2)^{\frac{5}{2}}}{(s-r)^3}ds =}
\nonumber \\ & \quad  &\qquad\qquad\quad
-4224r^{10} +12240r^8 -12768r^6 +5655r^4 -930r^2 +27 \ , \ \ |r|<1
\end{eqnarray}

\subsection{$I_{3}(U_n,3,r), \ n=0\cdots 6$}
\begin{equation}
\frac{1}{\pi}\int_{-1}^{1}\hypsng\frac{U_0(s)(1-s^2)^{\frac{5}{2}}}{(s-r)^3}ds
= -10r^3 + \frac{15}{2}r \ , \ \ |r|<1
\end{equation}

\begin{equation}
\frac{1}{\pi}\int_{-1}^{1}\hypsng\frac{U_1(s)(1-s^2)^{\frac{5}{2}}}{(s-r)^3}ds
= -30r^4 + 30r^2 - \frac{15}{4} \ , \ \ |r|<1
\end{equation}

\begin{equation}
\frac{1}{\pi}\int_{-1}^{1}\hypsng\frac{U_2(s)(1-s^2)^{\frac{5}{2}}}{(s-r)^3}ds
= -84r^5 + 110r^3 - 30r \ , \ \ |r|<1
\end{equation}

\begin{equation}
\frac{1}{\pi}\int_{-1}^{1}\hypsng\frac{U_3(s)(1-s^2)^{\frac{5}{2}}}{(s-r)^3}ds
= -224r^6 + 360r^4 - 150r^2 + 10 \ , \ \ |r|<1
\end{equation}

\begin{equation}
\frac{1}{\pi}\int_{-1}^{1}\hypsng\frac{U_4(s)(1-s^2)^{\frac{5}{2}}}{(s-r)^3}ds
= -576r^7 + 1092r^5 - 610r^3 + 90r \ , \ \ |r|<1
\end{equation}

\begin{equation}
\frac{1}{\pi}\int_{-1}^{1}\hypsng\frac{U_5(s)(1-s^2)^{\frac{5}{2}}}{(s-r)^3}ds
= -1440r^8 +3136r^6 -2190r^4 +510r^2 -20 \ , \ \ |r|<1
\end{equation}

\begin{equation}
\frac{1}{\pi}\int_{-1}^{1}\hypsng\frac{U_6(s)(1-s^2)^{\frac{5}{2}}}{(s-r)^3}ds
= -3520r^9 +8640r^7 -7224r^5 +2310r^3 -210r \ , \ \ |r|<1
\end{equation}


\subsection{$I_{4}(T_n,3,r), \ n=0\cdots 8$}
\begin{equation}
\frac{1}{\pi}\int_{-1}^{1}\hypsng\frac{T_0(s)(1-s^2)^{\frac{5}{2}}}{(s-r)^4}ds
= -10r^2 + \frac{5}{2} \ , \ \ |r|<1
\end{equation}

\begin{equation}
\frac{1}{\pi}\int_{-1}^{1}\hypsng\frac{T_1(s)(1-s^2)^{\frac{5}{2}}}{(s-r)^4}ds
= -20r^3 + 10r \ , \ \ |r|<1
\end{equation}

\begin{equation}
\frac{1}{\pi}\int_{-1}^{1}\hypsng\frac{T_2(s)(1-s^2)^{\frac{5}{2}}}{(s-r)^4}ds
= -70r^4 + 60r^2 - \frac{25}{4} \ , \ \ |r|<1
\end{equation}

\begin{equation}
\frac{1}{\pi}\int_{-1}^{1}\hypsng\frac{T_3(s)(1-s^2)^{\frac{5}{2}}}{(s-r)^4}ds
= -224r^5 + 260r^3 - 60r \ , \ \ |r|<1
\end{equation}

\begin{equation}
\frac{1}{\pi}\int_{-1}^{1}\hypsng\frac{T_4(s)(1-s^2)^{\frac{5}{2}}}{(s-r)^4}ds
= -672r^6 + 980r^4 - 360r^2 + 20 \ , \ \ |r|<1
\end{equation}

\begin{equation}
\frac{1}{\pi}\int_{-1}^{1}\hypsng\frac{T_5(s)(1-s^2)^{\frac{5}{2}}}{(s-r)^4}ds
= -1920r^7 +3360r^5 -1700r^3 +220r \ , \ \ |r|<1
\end{equation}

\begin{equation}
\frac{1}{\pi}\int_{-1}^{1}\hypsng\frac{T_6(s)(1-s^2)^{\frac{5}{2}}}{(s-r)^4}ds
= -5280r^8 +10752r^6 -6930r^4 +1460r^2 -50 \ , \ \ |r|<1
\end{equation}

\begin{equation}
\frac{1}{\pi}\int_{-1}^{1}\hypsng\frac{T_7(s)(1-s^2)^{\frac{5}{2}}}{(s-r)^4}ds
= -14080r^9 +32640r^7 -25536r^5 +7540r^3 -620r \ , \ \ |r|<1
\end{equation}

\begin{equation}
\frac{1}{\pi}\int_{-1}^{1}\hypsng\frac{T_8(s)(1-s^2)^{\frac{5}{2}}}{(s-r)^4}ds
= -36608r^{10} +95040r^8 -87360r^6 +33320r^4 -4560r^2 +104\ , \ \ |r|<1
\end{equation}

\subsection{$I_{4}(U_n,3,r), \ n=0\cdots 6$}
\begin{equation}
\frac{1}{\pi}\int_{-1}^{1}\hypsng\frac{U_0(s)(1-s^2)^{\frac{5}{2}}}{(s-r)^4}ds
= -10r^2 + \frac{5}{2} \ , \ \ |r|<1
\end{equation}

\begin{equation}
\frac{1}{\pi}\int_{-1}^{1}\hypsng\frac{U_1(s)(1-s^2)^{\frac{5}{2}}}{(s-r)^4}ds
= -40r^3 + 20r \ , \ \ |r|<1
\end{equation}

\begin{equation}
\frac{1}{\pi}\int_{-1}^{1}\hypsng\frac{U_2(s)(1-s^2)^{\frac{5}{2}}}{(s-r)^4}ds
= -140r^4 + 110r^2 - 10 \ , \ \ |r|<1
\end{equation}

\begin{equation}
\frac{1}{\pi}\int_{-1}^{1}\hypsng\frac{U_3(s)(1-s^2)^{\frac{5}{2}}}{(s-r)^4}ds
= -448r^5 + 480r^3 - 100r \ , \ \ |r|<1
\end{equation}

\begin{equation}
\frac{1}{\pi}\int_{-1}^{1}\hypsng\frac{U_4(s)(1-s^2)^{\frac{5}{2}}}{(s-r)^4}ds
= -1344r^6 + 1820r^4 - 610r^2 + 30 \ , \ \ |r|<1
\end{equation}

\begin{equation}
\frac{1}{\pi}\int_{-1}^{1}\hypsng\frac{U_5(s)(1-s^2)^{\frac{5}{2}}}{(s-r)^4}ds
= -3840r^7 +6272r^5 -2920r^3 +340r \ , \ \ |r|<1
\end{equation}

\begin{equation}
\frac{1}{\pi}\int_{-1}^{1}\hypsng\frac{U_6(s)(1-s^2)^{\frac{5}{2}}}{(s-r)^4}ds
= -10560r^8 +20160r^6 -12040r^4 +2310r^2 -70 \ , \ \ |r|<1
\end{equation}


\ 
\subsection{Others}
\begin{equation}
\int_{-1}^{1}(1-t^2)^{\frac{3}{2}}T_0(t) dt = \int_{-1}^{1}(1-t^2)^{\frac{3}{2}} dt
= \frac{3}{8}\pi
\end{equation}

\begin{equation}
\int_{-1}^{1}(1-t^2)^{\frac{3}{2}}T_2(t) dt = -\frac{\pi}{4}
\end{equation}

\begin{equation}
\int_{-1}^{1}(1-t^2)^{\frac{3}{2}}T_4(t) dt = \frac{\pi}{16}
\end{equation}

\begin{equation}
\int_{-1}^{1}(1-t^2)^{\frac{3}{2}}T_n(t) dt = 0 \ \ , \ \ \mbox{for
all}\ n ,\ \  \mbox{and}\ n\neq 0, 2 , 4 \ .
\end{equation}



\thanks{{\flushleft \Large{\textbf{Acknowledgements}}}\\
\\
We acknowledge the support from the USA
National Science Foundation (NSF) through grants CMS-9713798
(Mechanics \& Materials Program) and DMS-9600119 (Applied Mathematics Program).}

\bigskip
\bibliographystyle{plain}
\bibliography{abbrevs,fgm,gradient,snglr_num}


\end{document}